\documentclass[12pt]{amsart}
\usepackage{amssymb,amsmath,cite}

\newtheorem{proposition}{Proposition}[section]
\newtheorem{theorem}[proposition]{Theorem}

\newtheorem{lemma}[proposition]{Lemma}

\theoremstyle{definition}
\newtheorem{definition}[proposition]{Definition}

\theoremstyle{remark}
\newtheorem{remark}[proposition]{Remark}
\numberwithin{equation}{section}

\def\eps{\varepsilon}
\def\sig{\sigma}
\def\I {\mathbb{I}}
\def\R {\mathbb{R}}
\def\N {\mathbb{N}}
\def\C {\mathbb{C}}
\def\CC {{\mathcal C}}
\def\H {{\mathcal H}}
\def\D {{\mathcal D}}
\def\E {{\mathcal E}}
\def\F {{\mathcal F}}
\def\T {{\mathbb T}}
\def\X {{\mathcal X}}

\def\LL {{\mathbb L}}

\def\M {{\mathcal M}}
\def\Q {{\mathcal Q}}

\def\PP {{\mathbb P}}
\def\QQ {{\mathbb Q}}

\def\l {\langle}
\def\r {\rangle}
\def\sed {S_{\sigma,\tau,\eps}}
\def\soo {S_{0,0,0}}

\def\pep {{\boldsymbol{P_{\sig,\tau,\eps}}}}
\def \bra#1{\left\langle #1 \right\rangle}
\def\and{\qquad\text{and}\qquad}
\def\th{\vartheta}

\def \teta {\tilde \eta}

\def \au {\rm}
\def \ti {\it}
\def \jou {\rm}

\def \no#1#2#3 {{\bf #1} (#3), #2.}
\def \eds#1#2#3 {#1, #2, #3.}

\title[On a thermoviscoelastic plate with memory]{Asymptotic behavior
of a thermoviscoelastic \\ plate with memory effects}

\author[M.\ Grasselli]{Maurizio Grasselli}

\address{Dipartimento di Matematica ``F.Brioschi''
\newline\indent
Politecnico di Milano
\newline\indent
Via E.~Bonardi 9, 20133 Milano, Italy}
\email{maurizio.grasselli@polimi.it}

\author[J.E.\ Mu\~noz Rivera]{Jaime E. Mu\~noz Rivera}

\address{Laboratorio Nacional de Computa\c{c}ao Cientifica
\newline\indent
Av.~Get\'{u}lio Vargas 333, 25651-070-Petropolis, Brazil}
\email{rivera@lncc.br}

\author[M.\ Squassina]{Marco Squassina}
\address{Dipartimento di Informatica
\newline\indent
Universit\`a degli Studi di Verona
\newline\indent
C\'a Vignal 2, Strada Le Grazie 15, 37134 Verona, Italy}
\email{marco.squassina@univr.it}

\date{\today}

\thanks{The first author was partially supported by the
Italian PRIN Research Project 2006 {\it Problemi a frontiera libera,
transizioni di fase e modelli di isteresi}}

\thanks{The third author was partially supported by the
Italian PRIN Research Project 2007 {\it Metodi variazionali e topologici 
nello studio di fenomeni non lineari}}


\begin{document}

\begin{abstract}
We consider a coupled linear system describing a thermoviscoelastic plate
with hereditary effects. The system consists of a hyperbolic
integrodifferential equation, governing the temperature,
which is linearly coupled with the partial differential equation
ruling the evolution of the vertical
deflection, presenting a convolution term accounting for memory effects.
It is also assumed that the thermal power contains a memory term
characterized by a relaxation kernel.
We prove that the system is exponentially stable and we obtain a closeness estimate
between the system with memory effects and the corresponding memory-free limiting system,
as the kernels fade in a suitable sense.
\end{abstract}

\maketitle



\section{Introduction}
\noindent Let $\Omega$ be a bounded planar domain with smooth
boundary $\partial\Omega$. Suppose that $\Omega$ is occupied,
for all time $t$, by a thin homogeneous isotropic viscoelastic
plate. Denoting by $u$ its {\em vertical deflection} and by
$\vartheta$ the {\em temperature variation field}, we suppose
that the evolution of the pair $(u,\vartheta)$ is governed by
the following integrodifferential system
\begin{equation}
\label{orig-mob}
\begin{cases}
\displaystyle
u_{tt}+h(0)\Delta^2u+\int_0^\infty h'(s)\Delta^2u(t-s)ds+\Delta\vartheta=0, \\
\noalign{\vskip6pt}
\displaystyle
\vartheta_t+a(0)\vartheta+
\int_0^\infty a'(s)\vartheta(t-s)ds
-\displaystyle\int_0^\infty k(s)\Delta\vartheta(t-s)ds-\Delta u_t=0,
\end{cases}
\end{equation}
in $\Omega\times\R^+$, where $\R^+=(0,\infty)$. Here
$k:[0,\infty)\to\R^+$ and $h:[0,\infty)\to\R^+$ are smooth
decreasing convex functions which go to $0$ and to $h(\infty)>0$
at infinity, respectively. Instead, the memory kernel
$a:[0,\infty)\to\R^+$ is a smooth increasing concave function
with $a'$ vanishing at infinity. Moreover, all the other
physical constants have been set equal to $1$. Observe that, if
$k$ and $h-h(\infty)$ coincide with the Dirac mass $\delta_0$ at zero and
$a\equiv 0$, then, supposing $h(\infty)=1$, the above system formally collapses into the
linear model of thermoviscoelastic plate 
\begin{equation}
\label{orig-limit}
\begin{cases}
\displaystyle
u_{tt} +\Delta^2 u_t +\Delta(\Delta u+\vartheta)=0, \\
\noalign{\vskip5pt}
\vartheta_t-\Delta\vartheta-\Delta u_t=0.
\end{cases}
\end{equation}
We shall assume for simplicity that system \eqref{orig-mob} is endowed with Navier boundary conditions
\begin{align*}
u(t)=\Delta u(t)&=0 \qquad\text{on $\partial\Omega$,\, $t\geq 0$}, \\
\vartheta(t)&=0 \qquad\text{on $\partial\Omega$,\, $t\in\R$},
\end{align*}
and initial conditions
\begin{align*}
(u(0),u_t(0),\vartheta(0)) &=(u_0,u_1,\vartheta_0) \qquad\text{in $\Omega$}, \\
\vartheta(-s)&=\vartheta_0(s) \qquad\text{in $\Omega\times\R^+$}, \\
u(-s)&=u_0(s) \qquad\text{in $\Omega\times\R^+$},
\end{align*}
where $u_0,u_1,\vartheta_0:\Omega\to\R$ and $u_0,\vartheta_0:\Omega\times\R^+\to\R$
are assigned functions. The choice of these boundary conditions ({\em edge-free plate}) simplifies
the functional setup as well as some technical arguments with respect, e.g., to Neumann
boundary data ({\em clamped plate}). The results will be obtained via the so-called past history approach
(cf.\ \cite{Terreni} and references therein) which allows, under suitable assumptions, to express the solution by
a strongly continuous semigroup acting on an appropriate (extended) phase space
(cf.\ Theorem \ref{EU}).
\vskip2pt

System \eqref{orig-mob} with $h\equiv h(\infty)$ was considered 
and justified from the physical viewpoint in \cite{giopat},
while the viscoelastic case was treated in \cite{naso} (cf. also their references).
In \cite{giopat} the exponential stability was proved provided that $a(0)\not=0$, namely
not only the heat conduction law accounts for hereditary effects (see \cite{GuPi}),
but also the constitutive assumption for the thermal power contains a memory
term characterized by a nonzero relaxation kernel.
Instead, if one assumes that $h\equiv h(\infty)$ and $a\equiv 0$, then, for nonzero initial histories,
the system {\em fails} to be exponentially stable, no matter how fast the memory
kernel $k$ squeezes at infinity, provided that its growth around the origin is
suitably controlled (cf.\ \cite[Thm.~5.4]{GP-Riv}).
This confirms the conjecture that was formulated in \cite[Rem.~5.1]{giopat} and
also says that the presence of past history plays a discriminating role for
the stability of the thermoelastic system. It must be noticed that the
exponential stability was obtained in \cite{giopat} and in \cite{naso}
by exploiting some spectral analysis arguments, without detecting a precise decay rate. 
On the other hand, mainly in view of the asymptotic analysis that we wish to pursue
with respect to the behavior of the memory kernels involved in \eqref{orig-mob},
here we are interested in getting an explicit rate of decay.
By exploiting a technique first introduced in \cite{hazu,zua}, we detect the decay rate
by building up an ad hoc perturbation of the energy functional
which satisfy suitable differential inequalities (cf.\ Theorem \ref{edec}).
Concerning the case $h\not\equiv h(\infty)$ and $a\equiv 0$, in  Appendix we
shall consider a quite general class of
(abstract) thermoelastic systems with memory. We shall prove that
every trajectory squeezes to zero asymptotically (nonuniformly with respect to initial data).
Moreover, we shall exhibit some (weakly singular) memory kernels for
which the corresponding system, not including \eqref{orig-mob}, lacks of exponential stability.
\vskip2pt

The second main result of this paper is about the closeness between the solutions to
system \eqref{orig-mob} and the solutions to the system \eqref{orig-limit}.
The set of boundary and initial conditions is the same
but the ones for the past histories
of $\vartheta$ and $u$. Concerning the memory kernels $k$ and $h$, we proceed in the spirit
of \cite{cps1} (see also \cite{cps2,gnp,GS1}) by replacing them with the rescalings $k_\eps$ and $h_\sig$, defined by
\begin{equation*}
k_\eps(s)=\frac{1}{\eps}k\left(\frac{s}{\eps}\right), \qquad
h_\sig(s)=\frac{1}{\sig}\tilde h\left(\frac{s}{\sig}\right), \qquad \forall\, s\in\R^+,
\end{equation*}
where $\tilde h= h-h(\infty)$, while $\eps\in(0,1]$ and $\sig\in(0,1]$ are time relaxation parameters.
Notice that $k_\eps$ and $h_\sig$ approach the Dirac mass $\delta_0$ as $\eps$
and $\sig$ go to zero, in the sense of distribution. Moreover, on the basis
of physical motivations, concerning the parametrization
of the memory kernel $a$, 
we think of the (model) situation
$$
a_\tau(s)=\phi(\tau)+\psi(\tau)(1-e^{-\omega s}),
\qquad \forall s\in\R^+,
$$
$\phi,\psi:[0,1]\to\R^+$ being continuous functions with $\phi(0)=\psi(0)=0$
(see \eqref{nutauipsB}-\eqref{nuip}). Therefore, while the
kernels $k$ and $h$ undergo a singular perturbation procedure,
$a$ is parameterized just in order to be uniformly squeezing to zero as $\tau$ vanishes.
A suitable reformulation of system \eqref{orig-mob}, according to a well-established
procedure, with $k_\eps$ in place of $k$, $h_\eps$ in place of $h$,
and $a_\tau$ in place of $a$ is shown to generate a semigroup of contractions $\sed(t)$ on a certain
phase-space $\H_{\sig,\tau,\eps}^0$. Then, denoting by $\soo(t)$ the limiting
semigroup generated by system \eqref{orig-limit}, we establish an estimate
of the difference between two different trajectories, in terms of $\sig$, $\tau$ and $\eps$,
which holds on any bounded time interval. Basically, our estimate
says that the solutions to system \eqref{orig-mob} are arbitrarily close, in
the natural norm of $\H_{\sig,\tau,\eps}^0$, to the solutions of system \eqref{orig-limit},
provided that $\sig$, $\tau$ and $\eps$ are small enough and the initial data
are chosen inside a suitable regular bounded subset of the phase-space.
For the sake of generality we stress that, in addition to the singular limit
estimate in the norm of the base phase-space $\H_{\sig,\tau,\eps}^0$, we shall actually provide
the control with respect to the norms of the higher order phase-spaces
$\H_{\sig,\tau,\eps}^m$, $m\geq 0$ for suitably regular initial data (cf.\ Thm \ref{ThmGP1}).
Clearly, the limit process for $\sig$ and $\eps$ going to zero is {\em singular}, for the information on
the past histories of the temperature field $\vartheta$ and of the vertical deflection
$u$ get lost in the limit. As we shall see, the closeness control has to be
understood for time intervals which are bounded away from $0$.
In the particular case where we fix $\tau=0$ and we only take care of the limit process
with respect to $\sig$ and $\eps$, the result strengthens since
the estimate turns out to hold with constants which are independent of the time interval
size, so that the differences between any two trajectories can be controlled for any time $t>0$
(cf.\ Theorem \ref{ThmGP2}).

An interesting open problem is the analysis of the present model when $\tilde h$ is
approximated as $a$, that is, by a vanishing sequence of kernels. In fact, recalling that
there is no exponential decay when $a$ and $h'$ vanish (see \cite{GP-Riv}), there should
be a relation between the relaxation times $\eps$, $\sig$ and $\tau$ in order to preserve
the exponential stability when they approach zero.

\vskip4pt
The content of the paper is organized as follows.
\vskip4pt
\noindent
In \S\ref{setup} we introduce the notation and the basic tools, and we
formulate the problems in the proper functional setting.
In \S\ref{exps} we prove that, for every $\tau\neq 0$, the solutions
to \eqref{orig-mob} are exponentially decaying with a rate of decay
proportional to $\phi(\tau)$. 
In \S\ref{sing} we demonstrate the closeness estimate between
the strongly continuous semigroups associated with systems \eqref{orig-mob} and
\eqref{orig-limit} when the time rescaling parameters $\sig$, $\tau$ and $\eps$
tend to zero. Finally, in the Appendix, we deal with the pointwise decay and the
lack of exponential stability for an abstract class of thermoelastic
systems with memory.

\bigskip

\section{Preliminaries and well-posedness}
\label{setup}

\noindent
In this section we provide the proper functional
framework and the well-posedness result for problem \eqref{orig-mob}.

\subsection{The kernels parametrization}
We assume that $k:\R^+\to\R^+$ and $h:\R^+\to\R^+$ are smooth,
decreasing and summable functions satisfying, for the sake of
simplicity, the normalization conditions
$$
\int_0^\infty k(s)ds=
\int_0^\infty \tilde h(s)ds=1,\qquad
h(0)=2, \qquad k(0)=h(\infty)=1.
$$
Then, we set
$$
\mu(s)=-k'(s),\qquad \beta(s)=-h'(s),
\qquad \forall\, s\in\R^+,
$$
where $\mu$ and $\beta$ are supposed to satisfy
\begin{align}
\label{K0}
&\mu,\beta \in C^1(\R^+)\cap L^1(\R^+), \\
\label{K1}
& \mu(s)\ge 0, \quad\beta(s)\geq 0,
\qquad \forall\, s\in\R^+,\\
\noalign{\vskip4pt}
\label{K2}
& \mu'(s)\le 0, \quad\beta'(s)\leq 0,
\qquad \forall\, s\in\R^+,\\
\noalign{\vskip4pt}
\label{K3}
& \mu'(s) +\delta_1\mu(s) \leq 0,\qquad\forall\, s\in\R^+, \\
\noalign{\vskip5pt}
\label{K32}
& \beta'(s)+\delta_2\beta(s)\leq 0,\qquad\forall\, s\in\R^+.
\end{align}
for some $\delta_1>0$ and $\delta_2>0$. For any $\eps\in(0,1]$
and $\sig\in(0,1]$ we define the rescalings
\begin{equation}
\label{riscke}
\mu_\eps(s)=\frac{1}{\eps^2}\mu\left(\frac{s}{\eps}\right)=-k_\eps'(s), \qquad
\beta_\sigma(s)=\frac{1}{\sig^2}\beta\left(\frac{s}{\sig}\right)=-h_\sig'(s).
\end{equation}
Without loss of generality, we suppose
that, for $\eps\in(0,1]$ and $\sig\in(0,1]$, there holds
\begin{align}
\label{normalize1}
& \int_0^\infty \mu_\eps(s)ds=\textstyle{\frac{1}{\eps}},
\qquad
\displaystyle\int_0^\infty s\mu_\eps(s)ds=1,  \\
\label{normalize3}
& \int_0^\infty \beta_\sigma(s)ds=\textstyle{\frac{1}{\sig}},
\qquad
\displaystyle\int_0^\infty s\beta_\sigma(s)ds=1.
\end{align}
We assume that $a(s)=a_\tau(s)$, with $\tau\in[0,1]$, where $a_\tau:\R^+\to\R^+$
is a smooth concave function. We put $\nu_\tau(s)=-a''_\tau(s)$, where
$\nu_\tau \in C^1(\R^+)\cap L^1(\R^+)$ satisfies
\begin{align}
\label{nutauipsB}
\nu_\tau(s) \geq 0,\qquad
\nu'_\tau(s)\leq 0,\qquad\forall\,s\in\R^+, \\
\noalign{\vskip2pt}
\label{nutauips}
\nu'_\tau(s)+\delta_3\nu_\tau(s)\leq 0,\qquad\forall\, s\in\R^+,
\end{align}
for some $\delta_3>0$. Furthermore we assume that
the map $\{\tau\mapsto\nu_\tau\}$ is increasing
and there exist two functions $\phi,\psi\in C^0(\R^+)$ with $\phi\geq 0$, $\psi\geq 0$
and $\phi(0)=\psi(0)=0$, such that
\begin{equation}
\label{nuip}
a_\tau(0)=\phi(\tau),\qquad\forall\tau\in [0,1],
\qquad
\|\nu_\tau\|_{L^1(\R^+)}\leq \psi(\tau),\qquad\forall\tau\in [0,1].
\end{equation}

\subsection{The scale of phase-spaces}
Let $\Omega$ be a smooth bounded subset of $\R^2$.
The symbols $\|\cdot\|$ and $\langle\cdot,\cdot \rangle$ stand for the norm and the inner product on
$L^2(\Omega)$, respectively. We define the positive operator $A$ on $L^2(\Omega)$ by
$A=-\Delta$ with domain $\D(A)=H^1_0(\Omega)\cap H^2(\Omega)$,
and we introduce the scale of Hilbert spaces
$H^m=\D(A^{m/2})$, $m\in\R$, endowed with the inner products
$\langle u_1,u_2\rangle_ {H^m}=\langle A^{m/2}u_1,A^{m/2}u_2\rangle$.
We now consider the weighted Hilbert spaces
$$
\M_{\tau,\eps}^m=L^2_{\mu_\eps}(\R^+,H^{m+1})\cap L^2_{\nu_\tau}(\R^+,H^{m}),\qquad
\Q_\sigma^m=L^2_{\beta_\sigma}(\R^+,H^{m+1}),\qquad m\in\R,
$$
endowed, respectively, with the inner products
\begin{align*}
\l\eta_1,\eta_2\r_{\M_{\tau,\eps}^m}
&=\int_0^\infty\mu_\eps(s)\l A^{(1+m)/2}\eta_1(s),A^{(1+m)/2}\eta_2(s)\r ds \\
&+\int_0^\infty\nu_\tau(s)\l A^{m/2}\eta_1(s),A^{m/2}\eta_2(s)\r ds, \\
\noalign{\vskip2pt}
\l\xi_1,\xi_2\r_{\Q_\sig^m}
&=\int_0^\infty\beta_\sigma(s)\l A^{(1+m)/2}\xi_1(s),A^{(1+m)/2}\xi_2(s)\r ds,
\end{align*}
and we introduce the product spaces
\begin{align*}
\H^m_{\sigma,\tau,\eps} &=
\begin{cases}
H^{m+2}\times H^{m}\times H^m\times \M_{\tau,\eps}^m\times \Q_\sigma^{m+1},
& \text{if $\sig>0$ and $\tau>0$ or $\eps>0$},\\
\noalign{\vskip5pt}
H^{m+2}\times H^{m}\times H^m\times \M_{\tau,\eps}^m, & \text{if $\sig=0$ and $\tau>0$ or $\eps>0$},\\
\noalign{\vskip5pt}
H^{m+2}\times H^{m}\times H^m\times \Q_\sigma^{m+1}, & \text{if $\sig>0$ and $\tau=\eps=0$},\\
\noalign{\vskip5pt}
H^{m+2}\times H^{m}\times H^{m},  & \text{if $\sig=\tau=\eps=0$},
\end{cases}
\end{align*}
that will be normed by
\begin{equation*}
\|(u,u_t,\vartheta,\eta,\xi)\|_{\H^m_{\sigma,\tau,\eps}}^2=
\|u\|_{H^{m+2}}^2+\|u_t\|_{H^m}^2+\|\vartheta\|_{H^m}^2
+\|\eta\|_{\M^m_{\tau,\eps}}^2+\|\xi\|_{\Q^{m+1}_\sigma}^2.
\end{equation*}
In particular $\H^0_{\sigma,\tau,\eps}$ is the extended phase-space on which
we shall construct the dynamical system associated with~\eqref{orig-mob}.
Throughout the paper, when $\sig=\tau=\eps=0$, we shall agree to interpret the five entries vector
$z=(u,u_t,\vartheta,\eta,\xi)$ just as the triplet $(u,u_t,\vartheta)$.

\subsection{The problem setting}
In order to formulate the problem in a suitable
history space setting, let $T_{\tau,\eps}$ and $T_\sig$
be the linear operators on $\M^0_{\tau,\eps}$ and $\Q^1_\sig$
respectively, defined as
$$
T_{\tau,\eps}\eta=-\eta_s,\quad\eta\in{\D}(T_{\tau,\eps}),
\qquad
T_\sig\xi=-\xi_s,\quad\xi\in{\D}(T_\sig),
$$
where
\begin{align*}
\D(T_{\tau,\eps})&=\big\{\eta\in\M^0_{\tau,\eps}:\,\,
\eta_s\in \M^0_{\tau,\eps},\,\,\eta(0)=0\big\}, \\
\noalign{\vskip3pt}
\D(T_\sig) &=\big\{\xi\in\Q^1_\sig:\,\,
\xi_s\in \Q^1_\sig,\,\,\xi(0)=0\big\},
\end{align*}
and $\eta_s$ (resp.\ $\xi_s$) stands for the distributional derivative of
$\eta$ (resp.\ $\xi$) with respect to the internal variable $s$.
Notice that $T_{\tau,\eps}$ (resp.\ $T_\sigma$) is the infinitesimal generator of the right-translation
semigroup on $\M_{\tau,\eps}^0$ (resp.\ $\Q_\sig^1$). Moreover, on account
of \eqref{K2} and \eqref{nutauipsB},
\begin{align}
\label{etactrl}
\l T_{\tau,\eps}\eta,\eta\r_{\M_{\tau,\eps}^0}&= \frac{1}{2}\int_0^\infty \mu'_\eps(s)\Vert A^{1/2}\eta(s)\Vert^2 ds
+\frac{1}{2}\int_0^\infty \nu'_\tau(s)\Vert \eta(s)\Vert^2 ds\le 0, \\
\noalign{\vskip2pt}
\label{etactrl1}
\l T_\sig\xi,\xi\r_{\Q_\sig^1}&= \frac{1}{2}\int_0^\infty \beta'_\sig(s)\Vert A\xi(s)\Vert^2 ds \le 0,
\end{align}
for $\eta\in\D(T_{\tau,\eps})$ and $\xi\in\D(T_\sig)$.
Following the well-established past history approach (see, e.g.\ \cite{Terreni}),
we introduce the so-called past histories of $\vartheta$ and $u$,
\begin{equation*}
\eta^t(s)=\int_{0}^s \vartheta(t-y) dy,\qquad
\xi^t(s)=u(t)-u(t-s),
\qquad (s,t)\in \R^+\times\R^+.
\end{equation*}
Differentiating these variables leads to
further equations ruling the evolution of $\eta$ and $\xi$
\begin{equation*}
\eta^t_t=-\eta^t_s+\vartheta(t),\qquad
\xi^t_t=-\xi^t_s+u_t(t),
\qquad t\in\R^+.
\end{equation*}
We are now in the right position to introduce the formulation of the problems.
On account of the normalization conditions and of the notation previously introduced,
for any $\sig,\tau,\eps\in[0,1]$, given $(u_0,u_1,\vartheta_0,\eta_0,\xi_0)$ in $\H^0_{\sigma,\tau,\eps}$,
find $(u,u_t,\vartheta,\eta,\xi)\in C([0,\infty),\H^0_{\sig,\tau,\eps})$ solution to
\begin{equation}
\tag{${\mathcal P}_{\sigma,\tau,\eps}$}
\label{eq-mobile}
\begin{cases}
\displaystyle
u_{tt}+\int_0^\infty \beta_\sigma(s)A^2\xi(s)ds+A(Au-\vartheta)=0, \\
\noalign{\vskip5pt}
\displaystyle
\vartheta_t+\phi(\tau)\vartheta
+\int_0^\infty \nu_\tau(s)\eta(s)ds
+\int_0^\infty \mu_\eps(s)A\eta(s)ds+Au_t=0, \\
\noalign{\vskip6pt}
\eta_t=T_{\tau,\eps}\eta+\vartheta, \\
\noalign{\vskip5pt}
\xi_t=T_\sigma\xi+u_t,
\end{cases}
\end{equation}
for $t\in\R^+$, with initial condition
$(u(0),u_t(0),\vartheta(0),\eta^0,\xi^0)=(u_0,u_1,\vartheta_0,\eta_0,\xi_0)$.
Similarly, we introduce the limiting problem (formally corresponding
to the case $\sigma=\tau=\eps=0$).
\noindent
Given $(u_0,u_1,\vartheta_0)\in \H^0_{0,0,0}$, find $(u,u_t,\vartheta)\in C([0,\infty),
\H^0_{0,0,0})$ solution to
\begin{equation}
\tag{${\mathcal P}_{0,0,0}$}
\label{eq-limite}
\begin{cases}
\displaystyle
u_{tt}+A^2u_t+A(Au-\vartheta)=0, \\
\noalign{\vskip5pt}
\vartheta_t+A\vartheta+Au_t=0,
\end{cases}
\end{equation}
for $t\in\R^+$, which fulfills the initial conditions
$(u(0),u_t(0),\vartheta(0))=(u_0,u_1,\vartheta_0)$.
The above problems are abstract reformulation of the initial and
boundary value problems associated with \eqref{orig-mob} and \eqref{orig-limit}.

\subsection{Well-posedness}
System ${\mathcal P}_{\sig,\tau,\eps}$
allows us to provide a description of the solutions in terms of a
strongly continuous semigroup of operators on $\H^0_{\sig,\tau,\eps}$.
Indeed, setting 
$$
\zeta(t)=(u(t),v(t),\vartheta(t),\eta^t,\xi^t)^{\top},
$$
the problem rewrites as
\begin{equation*}
\frac{d}{dt}\zeta={\mathcal L}\zeta, \qquad
\zeta(0)=\zeta_0,
\end{equation*}
where ${\mathcal L}$ is the linear operator defined by
\begin{equation}
\label{cauchy}
{\mathcal L}
\begin{pmatrix}
u \\
v \\
\vartheta \\
\eta \\
\xi
\end{pmatrix}
=
\begin{pmatrix}
v \\
-\int_0^\infty \beta_\sigma(s)A^2\xi(s)ds-A(Au-\vartheta) \\
\noalign{\vskip4pt}
-\phi(\tau)\vartheta-\int_0^\infty \nu_\tau(s)\eta(s)ds
-\int_0^\infty\mu_\eps(s)A\eta(s)ds-Av \\
\noalign{\vskip4pt}
\vartheta +T_{\tau,\eps}\eta \\
\noalign{\vskip4pt}
v +T_\sig\xi
\end{pmatrix}
\end{equation}
with domain
$$
{\mathcal D}({\mathcal L})=
\left\{z\in\H^0_{\sig,\tau,\eps}
\left|
\begin{aligned}
& Au-\vartheta\in H^2 \\
& v\in H^2,\,\,\, \vartheta \in H^1  \\
& \textstyle\int_0^\infty\mu_\eps(s)A\eta(s)ds\in H^0 \\
& \textstyle\int_0^\infty\nu_\tau(s)\eta(s)ds\in H^0 \\
& \textstyle\int_0^\infty\beta_\sig(s)A^2\xi(s)ds\in H^0 \\
& \eta\in{\D}(T_{\tau,\eps}),\,\,\, \xi\in{\D}(T_\sig)
\end{aligned}
\right\}.
\right.
$$
\vskip4pt
\noindent
By virtue of \eqref{etactrl} and \eqref{etactrl1}, it is readily seen that ${\mathcal L}$ is a
dissipative operator. We will tacitly extend the definition of $\sed(t)$ to the case
$\sig=\tau=\eps=0$ which is well known. Of course, in this case the solution
semigroup is a three-component vector only. We now assume that \eqref{K0}-\eqref{K1},
\eqref{K3}-\eqref{K32} and \eqref{nutauipsB}-\eqref{nuip} hold true. If $\sig>0$, $\tau>0$ and $\eps>0$, 
following the proof of  \cite[Thm.~2.1]{giopat}, we obtain 

\begin{theorem}
\label{EU}
System ${\mathcal P}_{\sig,\tau,\eps}$ defines a $C_0$-semigroup $\sed(t)$ of
contractions on $\H^0_{\sig,\tau,\eps}$.
\end{theorem}

\bigskip

\section{Exponential stability of $\sed(t)$}
\label{exps}

\noindent
In this section we prove that, for any $\tau\in[0,1]$, the semigroup $\sed(t)$
is exponentially stable on $\H^0_{\sig,\tau,\eps}$, admitting
a decay rate proportional to $\phi(\tau)$ when $\tau>0$.
In this case, the exponential stability is actually already known
from \cite{giopat} in the elastic case with a nonvanishing kernel $a$ (recall that
if $a$ vanishes the exponential stability fails as shown in \cite{GP-Riv}).
However, this result was proven via spectral analysis arguments,
without detecting a precise decay rate.
Here, we exploit a technique first introduced in \cite{hazu,zua}, namely, and
we obtain the decay estimate for a suitably defined perturbation of the energy functional
$\E:\R^+\to\R$, defined by $\E(t)=\|\sed(t)\|^2_{\H^0_{\sig,\tau,\eps}}$. We can thus
provide a decay rate which shows the role played by the kernels $a$ and $h$.
\vskip4pt
The main result of this section is the following

\begin{theorem}
\label{edec}
Assume that \eqref{K0}-\eqref{K32} and \eqref{nutauipsB}-\eqref{nutauips}
hold. Then there exist $\Theta>0$, $d_0>0$ and $\varsigma>0$,
independent of $\sig$, $\tau$ and $\eps$, such that for any $\tau\in[0,1]$
\begin{equation}
\label{explicitdecay}
\E(t)\leq \varsigma \E(0)e^{-(\phi(\tau)+d_0)\Theta t},\qquad\forall\, t\geq 0.
\end{equation}
\end{theorem}

\begin{remark}
A careful analysis of the proof of Theorem~\ref{edec} shows that
the constant $d_0$ and $\Theta$ can be explicitly calculated. In particular, $d_0$ accounts for
the viscoelastic effects.
\end{remark}

\begin{proof}
Let $\tau\in[0,1]$ and let $0<\rho_\flat<\rho_\sharp<1$ to be chosen later. Then, for all $t\geq 0$,
consider the following perturbation $\F_1$ of the energy functional $\E$
\begin{equation*}
\F_1(t)=\E(t)+\rho_\flat\Theta_\flat(t)+\rho_\sharp\Theta_\sharp(t),
\qquad
\Theta_\flat(t)=\l u_t(t),u(t)\r,
\quad
\Theta_\sharp(t)=-\sigma\l u_t(t),\xi^t\r_{\Q_\sigma^{-1}}.
\end{equation*}
We denote by $C$ a generic positive constant independent of $\rho_\flat,\rho_\sharp$ and $\sigma,\tau,\eps$
which may vary from line to line within the same formula. Observe that, by \eqref{normalize3}, there holds
$$
|\Theta_\flat(t)|+|\Theta_\sharp(t)|\leq C\big(\|Au(t)\|^2
+\|u_t(t)\|^2+\|\xi^t\|_{\Q_\sigma^1}^2\big)
\leq C\E(t).
$$
Therefore, up to choosing $\rho_\flat$ and $\rho_\sharp$ sufficiently small, we have
$\frac{1}{2}\F_1(t)\leq \E(t)\leq 2 \F_1(t)$,
so that $\E$ and $\F_1$ turn out to be equivalent for what
concerns the energy decay estimate.
Let us now multiply the first equation of
${\mathcal P}_{\sig,\tau,\eps}$ by $u_t$ in $H^0$, the second
by $\vartheta$ in $H^0$, the third by $\eta$ in $\M_{\tau,\eps}^0$,
the fourth by $\xi$ in $\Q_\sig^1$ and
add the resulting identities. This yields
\begin{align*}
\frac{d}{dt}\E(t) & \leq -\textstyle{\frac{\delta_1}{\eps}}\|\eta^t\|_{L^2_{\mu_\eps}(\R^+,H^1)}^2
-\delta_3\|\eta^t\|_{L^2_{\nu_\tau}(\R^+,H^0)}^2
-\textstyle{\frac{\delta_2}{2\sigma}}\|\xi^t\|_{\Q_\sig^1}^2 \\
\noalign{\vskip3pt}
&+\textstyle{\frac{1}{2}}\displaystyle\int_0^\infty \beta_\sigma'(s)\|A\xi^t(s)\|^2ds
-2\phi(\tau)\|\vartheta(t)\|^2,
\end{align*}
by virtue of inequalities \eqref{K3}-\eqref{K32}, \eqref{nutauips}, \eqref{etactrl}-\eqref{etactrl1}
and integration by parts. Besides, by direct computation, we get
\begin{align*}
\frac{d}{dt}\Theta_\flat(t) &=\|u_t(t)\|^2-\|Au(t)\|^2+\l \vartheta(t),Au(t)\r
-\l\xi^t,u(t) \r_{\Q_\sig^1}, \\
\noalign{\vskip2pt}
\frac{d}{dt}\Theta_\sharp(t)&=-\sigma\l u_{tt}(t),\xi^t\r_{\Q_\sig^{-1}}
-\sigma\l u_t(t), T_\sig\xi^t\r_{\Q_\sig^{-1}}
-\|u_t(t)\|^2,
\end{align*}
where, in the last identity, we have used formula \eqref{normalize3} once again.
Then, on account of the obtained formulas for the derivatives of
$\E$, $\Theta_\flat$ and $\Theta_\sharp$, we deduce
\begin{align*}
\frac{d}{dt}\F_1(t) &\leq -{\textstyle
\min\{\delta_1/\eps,\delta_3\}}\|\eta^t\|_{\M_{\tau,\eps}^0}^2
-{\textstyle\frac{\delta_2}{2\sigma}}\|\xi^t\|_{\Q_\sig^1}^2
+{\textstyle\frac{1}{2}}\int_0^\infty \beta_\sigma'(s)\|A\xi^t(s)\|^2ds \\
\noalign{\vskip6pt}
&-2\phi(\tau)\|\vartheta(t)\|+\rho_\flat\|u_t(t)\|^2-\rho_\flat\|Au(t)\|^2+\rho_\flat\l \vartheta(t),Au(t)\r \\
\noalign{\vskip8pt}
&-\rho_\flat\l\xi^t,u(t) \r_{\Q_\sigma^1}
-\rho_\sharp\sigma\l u_{tt}(t),\xi^t\r_{\Q_\sig^{-1}}
-\rho_\sharp\sigma\l u_t(t), T_\sig\xi^t\r_{\Q_\sig^{-1}}
-\rho_\sharp\|u_t(t)\|^2.
\end{align*}
Therefore, we get
\begin{align*}
\frac{d}{dt}\F_1(t) &\leq -\rho_\flat\|Au(t)\|^2
-(\rho_\sharp-\rho_\flat)\|u_t(t)\|^2
-2\phi(\tau)\|\vartheta(t)\|^2 \\
\noalign{\vskip3pt}
&-\min\{\delta_1/\eps,\delta_3\}\|\eta^t\|_{\M_{\tau,\eps}^0}^2
-\textstyle{\frac{\delta_2}{2\sigma}}\|\xi^t\|_{\Q_\sig^1}^2
+{\textstyle\frac{1}{2}}\displaystyle\int_0^\infty \beta_\sigma'(s)\|A\xi^t(s)\|^2ds+{\mathcal J}(t),
\end{align*}
where we have set
\begin{equation*}
{\mathcal J}(t)=\rho_\flat\l \vartheta(t),Au(t)\r
-\rho_\flat\l u(t),\xi^t \r_{\Q_\sigma^1}
+\rho_\sharp\sigma\l u_t(t), \xi^t_s\r_{\Q_\sig^{-1}}
-\rho_\sharp\sigma\l u_{tt}(t),\xi^t\r_{\Q_\sig^{-1}}.
\end{equation*}
Notice that, we have
\begin{align*}
\l \vartheta(t),Au(t)\r &\leq \|\vartheta(t)\|^2+\textstyle{\frac{1}{4}}\|Au(t)\|^2, \\
\noalign{\vskip6pt}
-\l u(t),\xi^t \r_{\Q_\sigma^1} & \leq C\rho_\flat\|Au(t)\|^2+
\textstyle{\frac{\delta_2}{8\rho_\flat\sig}} \|\xi^t\|_{\Q_\sig^1}^2, \\
\noalign{\vskip4pt}
\sigma\l u_t(t), \xi^t_s\r_{\Q_\sig^{-1}} &\leq
C\rho_\sharp\|u_t(t)\|^2-\textstyle{\frac{1}{2\rho_\sharp}}
\displaystyle\int_0^\infty \beta_\sigma'(s)\|A\xi^t(s)\|^2ds, \\
-\sigma\l u_{tt}(t),\xi^t\r_{\Q_\sig^{-1}}&=\sigma\l u(t),\xi^t\r_{\Q_\sig^1}
-\sigma\l \vartheta(t),\xi^t\r_{\Q_\sig^{0}}+
\sigma\Big\|\int_0^\infty \beta_\sigma(s)A\xi^t(s)ds\Big\|^2 \\
\noalign{\vskip2pt}
& \leq C\rho_\sharp\|Au(t)\|^2+C\rho_\sharp\|\vartheta(t)\|^2+
\textstyle{\frac{\delta_2}{8\rho_\sharp\sigma}}\|\xi^t\|_{\Q_\sig^1}^2
+C\|\xi^t\|_{\Q_\sig^1}^2.
\end{align*}
By the above inequalities, it follows
\begin{align*}
{\mathcal J}(t) & \leq \big(\textstyle{\frac{\rho_\flat}{4}+C\rho_\flat^2+C\rho_\sharp^2}\big)\|Au(t)\|^2+C\rho_\sharp^2\|u_t(t)\|^2
+\big(\rho_\flat+C\rho_\sharp^2\big)\|\vartheta(t)\|^2 \\
\noalign{\vskip4pt}
&+\big(\textstyle{\frac{\delta_2}{4\sig}+C\rho_\sharp}\big)\|\xi^t\|_{\Q_{\sig}^1}^2
-\textstyle{\frac{1}{2}}\displaystyle\int_0^\infty \beta_\sigma'(s)\|A\xi^t(s)\|^2ds
\end{align*}
Therefore, we conclude that
\begin{align*}
\frac{d}{dt}\F_1(t)&+\big(\textstyle{\frac{3\rho_\flat}{4}-C\rho_\flat^2-C\rho_\sharp^2}\big)\|Au(t)\|^2 \\
\noalign{\vskip4pt}
&+(\rho_\sharp-\rho_\flat-C\rho_\sharp^2)\|u_t(t)\|^2
+\big(2\phi(\tau)-\rho_\flat-C\rho_\sharp^2\big)\|\vartheta(t)\|^2  \\
\noalign{\vskip5pt}
&+\textstyle{\min\{\delta_1/\eps,\delta_3\}}\|\eta^t\|_{\M_{\eps}^0}^2
+\textstyle{\frac{\delta_2-C\rho_\sharp\sig}{4\sig}}\|\xi^t\|_{\Q_\sig^1}^2\leq 0.
\end{align*}
Choosing $\rho_\flat=\bar\rho_\flat\phi(\tau)$ and
$\rho_\sharp=\bar\rho_\sharp\phi(\tau)$, where the positive constants
$\bar\rho_\flat$ and $\bar\rho_\sharp$ are independent
of $\sig,\tau$ and $\eps$, we obtain 
\begin{align*}
\frac{d}{dt}\F_1(t)&+\phi(\tau)(\textstyle{\frac{3\bar\rho_\flat}{4}-C\bar\rho_\flat^2-C\bar\rho_\sharp^2})\|Au(t)\|^2 \\
\noalign{\vskip4pt}
&+\phi(\tau)(\bar\rho_\sharp-\bar\rho_\flat-C\bar\rho_\sharp^2)\|u_t(t)\|^2
+\phi(\tau)(2-\bar\rho_\flat-C\bar\rho_\sharp^2)\|\vartheta(t)\|^2 \\
\noalign{\vskip5pt}
&+\phi(\tau)\textstyle{\min\{\delta_1,\delta_3\}}C\|\eta^t\|_{\M_{\tau,\eps}^0}^2
+\phi(\tau)\textstyle{\frac{\delta_2-C\bar\rho_\sharp}{4}}C\|\xi^t\|_{\Q_\sig^1}^2\leq 0.
\end{align*}
Then, fixing $\bar\rho_\flat$ and $\bar\rho_\sharp$ so small that
\begin{equation*}
\Lambda=\min\Big\{\textstyle{\frac{3\bar\rho_\flat}{4}-C\bar\rho_\flat^2-C\bar\rho_\sharp^2},
\bar\rho_\sharp-\bar\rho_\flat-C\bar\rho_\sharp^2,
2-\bar\rho_\flat-C\bar\rho_\sharp^2,
\textstyle{\min\{\delta_1,\delta_3\}C},
\textstyle{\frac{\delta_2-C\bar\rho_\sharp}{4}C}\Big\}>0,
\end{equation*}
and $\E$ controls and it is controlled
by $\F_1$, it follows that
\begin{equation}
\label{PRIM}
\frac{d}{dt}\F_1(t)+{\textstyle\frac{\Lambda}{2}\phi(\tau)}
 \F_1(t)\leq 0,\qquad\forall\, t\geq 0.
\end{equation}
By arguing as above, we have
\begin{equation}
\label{eq7}
\frac{d}{dt}\Theta_\sharp(t)\leq  -\frac{1}{2} \|u_t\|^2 +\frac{1}{32}\|Au\|^2+\frac{1}{16}\|\th\|^2
+\frac{C}{\sigma}\int_0^\infty\beta_\sigma(s)\|A\xi^t(s)\|^2ds.
\end{equation}
Setting now
$$
K(t)=-\eps\l\th(t),\eta^t\r_{\M_{0,\eps}^{-1}},
$$ 
multiplying the second equation of \eqref{eq-mobile}
by $\int_0^\infty\mu_\eps(s)\eta^t(s)ds$ and recalling \eqref{normalize1}, we get
\begin{align*}
\frac{d}{dt}K(t)&=\eps\Big\|\int^{\infty}_0 \mu_\eps(s)A^{1/2}\eta^t(s)ds\Big\|^2
+\frac{d}{dt} \eps\Big\langle Au,\int_0^\infty\mu_\eps(s)\eta^t(s)ds\Big\rangle \\
\noalign{\vskip2pt}
&\qquad   -\eps\Big\langle Au,\int_0^\infty\mu_\eps(s)\eta^t_t(s)ds\Big\rangle
+\eps\Big\langle\int_0^\infty\nu_\tau(s)\eta^t(s)ds,\int_0^\infty\mu_\eps(s)\eta^t(s)ds\Big\rangle       \\
\noalign{\vskip2pt}
&\qquad -\eps\Big\langle\th,\int_0^\infty\mu'_\eps(s)\eta^t(s)ds\Big\rangle-\|\th\|^2.
\end{align*}
Therefore, setting
$$
K_2(t)=K(t)- \eps\Big\langle Au,\int_0^\infty\mu_\eps(s)\eta^t(s)ds\Big\rangle,
$$
we obtain
\begin{align*}
\frac{d}{dt}K_2(t) &=\eps\Big\|\int^{\infty}_0 \mu_\eps(s)A^{1/2}\eta^t(s)ds\Big\|^2-\eps\Big\langle Au,\int_0^\infty\mu'_\eps(s)\eta^t(s)ds
\Big\rangle-\langle Au,\th \rangle \\
& +\eps\Big\langle\int_0^\infty\nu_\tau(s)\eta^t(s)ds,\int_0^\infty\mu_\eps(s)\eta^t(s)ds\Big\rangle  
-\eps\Big\langle \th,\int_0^\infty\mu'_\eps(s)\eta^t(s)ds\Big\rangle-\|\th\|^2.
\end{align*}
Notice that we get
\begin{align*}
-\eps\Big\langle \th,\int_0^\infty\mu'_\eps(s)\eta^t(s)ds\Big\rangle &\leq
\frac{1}{2}\|\th(t)\|^2+\frac{1}{2}\left[ \int_0^\infty\frac{-\eps\mu'_\eps(s)}{\mu_{\eps}^{1/2}(s)}\mu_{\eps}^{1/2}(s)\|\eta^t(s)\|ds  \right]^2 \\
&\leq \frac{1}{2}\|\th(t)\|^2+\frac{1}{2}\int_0^\infty\frac{\eps^2(\mu'_\eps(s))^2}{\mu_{\eps}(s)}ds\int_0^\infty\mu_{\eps}(s)\|\eta^t(s)\|^2ds \\
&\leq \frac{1}{2}\|\th(t)\|^2+\frac{C}{\eps}\int_0^\infty\mu_{\eps}(s)\|\eta^t(s)\|^2ds.
\end{align*}
Moreover,
\begin{gather*}
\eps\Big\langle\int_0^\infty\nu_\tau(s)\eta^t(s)ds,\int_0^\infty\mu_\eps(s)\eta^t(s)ds\Big\rangle \leq
\frac{\eps}{2}\left\| \int_0^\infty\nu_\tau(s)\eta^t(s)ds  \right\|^2+\frac{\eps}{2}\left\|\int_0^\infty\mu_\eps(s)\eta^t(s)ds  \right\|^2 \\
  \leq \psi(\tau)\int^{\infty}_0 \nu_\tau(s)\|\eta^t(s)\|^2ds+\int^{\infty}_0 \mu_\eps(s)\|A^{1/2}\eta^t(s)\|^2ds.
\end{gather*}
Hence we deduce the following inequality
\begin{align}
\label{eq8}
\frac{d}{dt}K_2(t) &\leq -\frac{1}{2}\|\th\|^2+\frac{C}{\eps}\int^{\infty}_0 \mu_\eps(s)\|A^{1/2}\eta^t(s)\|^2ds \\
&+\psi(\tau)\int^{\infty}_0 \nu_\tau(s)\|\eta^t(s)\|^2ds +\frac{1}{4}\|Au\|^2-\langle Au,\th \rangle. \notag
\end{align}
By multiplying the first equation of \eqref{eq-mobile} by $u$ we get
\begin{equation*}
\frac{d}{dt}
\left(u_{t},u\right)\leq \|u_{t}\|^2+\frac{1}{2\sigma}\int_0^\infty\beta_\sigma(s)\|A\xi^t(s)\|^2ds
-\frac 12\|Au\|^2+\left\langle \th,Au\right\rangle.
\end{equation*}
Whence, we deduce that
\begin{align*}
\frac{d}{dt}\big[K_2(t)+\left(u_{t},u\right)\big] &\leq
-\frac{1}{4}\|Au\|^2-\frac{1}{2}\|\th\|^2+\|u_{t}\|^2+C\int^{\infty}_0 \nu_\tau(s)\|\eta^t(s)\|^2ds   \\
\noalign{\vskip4pt}
&+\frac{1}{2\sigma}\int_0^\infty\beta_\sigma(s)\|A\xi^t(s)\|^2ds
+\frac{C}{\eps}\int^{\infty}_0 \mu_\eps(s)\|A^{1/2}\eta^t(s)\|^2ds.
\end{align*}
Using now inequality \eqref{eq7} and setting
$K_3(t)=4\Theta_\sharp(t)+K_2(t)+\left(u_{t},u\right)$, we get
\begin{align*}
\frac{d}{dt}K_3(t)
&\leq
-\frac {1}{8}\|Au\|^2-\frac{1}{4} \|\th\|^2-\|u_{t}\|^2+C\int^{\infty}_0 \nu_\tau(s)\|\eta^t(s)\|^2ds  \\
\noalign{\vskip4pt}
&+\frac{C}{\sigma}\int_0^\infty\beta_\sigma(s)\|A\xi^t(s)\|^2ds
+\frac{C}{\eps}\int^{\infty}_0 \mu_\eps(s)\|A^{1/2}\eta^t(s)\|^2ds.
\end{align*}
Since, as can be readily checked, it holds
$$
\frac{d}{dt}\E(t) \leq -\textstyle{\frac{\delta_1}{\eps}}\|\eta^t\|_{L^2_{\mu_\eps}(\R^+,H^1)}^2
-\delta_3\|\eta^t\|_{L^2_{\nu_\tau}(\R^+,H^0)}^2-\textstyle{\frac{\delta_2}{\sigma}}\|\xi^t\|_{\Q_\sig^1}^2,
$$
setting $\F_2(t)=N\mathcal{E}(t)+K_3(t)$ with $N$ sufficiently large and independent of $\eps$ and $\sigma$,
it is readily seen that $\F_2$ controls and it is controlled by the energy and
\begin{equation}
\label{SEC}
\frac{d}{dt}\F_2(t)+d_0\F_2(t)\leq 0,\qquad\forall t\geq 0,
\end{equation}
for some positive constant $d_0$ independent of $\eps$ and $\sigma$.
Therefore,  by combining inequalities \eqref{PRIM} and \eqref{SEC} and setting
$\F=\F_1+\F_2$ it follows that $\F$ is equivalent to the energy and satisfies
\begin{equation*}
\frac{d}{dt}\F(t)+C(\phi(\tau)+d_0)\F(t)\leq 0,\qquad\forall\, t\geq 0.
\end{equation*}
By the Gronwall Lemma we obtain the desired inequality \eqref{explicitdecay}.
\end{proof}

\begin{remark}
We know that
$$
\lim_{t\to\infty} \|S_{\sig,0,\eps}(t)z\|_{\H_{\sig,0,\eps}}=0,
\qquad\forall \, z\in\H_{\sig,0,\eps},
$$
provided that $\mu$ satisfies a mild summability condition (see
Theorem \ref{dipp} in Appendix). If, in addition, we assume
that $\beta\equiv 0$ and the memory kernel $\mu$
does not grow too rapidly around the origin, i.e.\ $\sqrt{s}\,\mu(s)\to 0$
for $s\to 0$, then the first order energy fails to vanish
through an exponential law of decay (cf.\ \cite[Thm.~5.4]{GP-Riv}
as well as Theorem \ref{ThmMAIN0} in Appendix for a more general situation).
For the case general $\beta\neq 0$, we refer the reader to the
Appendix for a discussion on the lack of exponential stability
for a class of abstract linear thermoelastic systems with (possibly) fractional
operator powers (cf.\ Theorem \ref{ThmMAIN}).
\end{remark}

\bigskip

\section{Closeness between $\sed(t)$ and $\soo(t)$}
\label{sing}

\noindent
The aim of this section is to establish, following a pattern recently initiated in \cite{cps1},
a precise quantitative estimate of the closeness between the {\em non analytic semigroup} $\sed(t)$ and 
the {\em analytic semigroup} $\soo(t)$ (see \cite{LaTr}) in the
norm of any extended phase-space $\H_{\sig,\tau,\eps}^m$, for $m\geq 0$,
as the parameters $\sig$, $\tau$ and $\eps$ converge to zero,
provided that the initial data are chosen inside a suitable
regular bounded subset of $\H_{\sig,\tau,\eps}^m$ (see also \cite{cps2,gnp}).
In \cite{GS1} a similar analysis was carried on in the case $\sig=\tau=0$,
for a plate model which accounts for the rotational inertia term $-\Delta u_{tt}$
in the equation ruling the vertical deflection.
We point out that, along the convergence process, $\mu_\eps$ and $\beta_\sig$
behave in a singular fashion since $\|\mu_\eps\|_{L^1(\R^+)}\to\infty$ and
$\|\beta_\sig\|_{L^1(\R^+)}\to\infty$ for $\eps$ and $\sig$ going to zero,
whereas the kernel $\nu_\tau$ satisfies $\|\nu_\tau\|_{L^1(\R^+)}\to 0$ as $\tau$ vanishes.

\subsection{Discussion of the results}
\label{introsing}
Throughout the section we will assume that, whenever $\sig>0$, $\tau>0$ and $\eps>0$,
conditions \eqref{K0}-\eqref{K1}, \eqref{K3}-\eqref{K32}
and \eqref{nutauipsB}-\eqref{nuip} hold true.
In order to perform a comparison between the five component
semigroup $\sed(t)$ and the three component (for
$\sig=\tau=\eps=0$) limiting semigroup $\soo(t)$,
we need to introduce, for any $m\geq 0$, the following lifting and projection maps
\begin{align*}
\LL_{\sig,\tau,\eps}&: \H_{0,0,0}^m\to\H_{\sig,\tau,\eps}^m, \\
\noalign{\vskip3pt}
\PP&: \H_{\sig,\tau,\eps}^m\to\H_{0,0,0}^m, \\
\noalign{\vskip3pt}
\QQ_{\tau,\eps}&: \H_{\sig,\tau,\eps}^m\to\M_{\tau,\eps}^m, \\
\noalign{\vskip3pt}
\QQ_{\sigma}&: \H_{\sig,\tau,\eps}^m\to\Q_\sigma^{m+1},
\end{align*}
defined, respectively, by
$$
\LL_{\sig,\tau,\eps} (u,u_t,\vartheta)=
\begin{cases}
(u,u_t,\vartheta,0,0),&\text{if $\sig>0$ and $\tau>0$ or $\eps>0$},\\
\noalign{\vskip1pt}
(u,u_t,\vartheta,0),&\text{if $\sig=0$ or $\tau=\eps=0$},\\
\noalign{\vskip1pt}
(u,u_t,\vartheta),&\text{if $\sig=\tau=\eps=0$},
\end{cases}
$$
and by
\begin{align*}
\PP(u,u_t,\vartheta,\eta,\xi)& =(u,u_t,\vartheta), \\
\noalign{\vskip3pt}
\QQ_{\tau,\eps}(u,u_t,\vartheta,\eta,\xi)& =\eta, \\
\noalign{\vskip3pt}
\QQ_{\sigma}(u,u_t,\vartheta,\eta,\xi)& =\xi.
\end{align*}
In the case $\tau>0$, if $z$ denotes the initial data, taken
inside any bounded subset of $\H_{\sig,\tau,\eps}^{2m+4}$, we will prove
the convergence of $\sed (t)z$ towards $\LL_{\sig,\tau,\eps}\soo(t)
\PP z$ in the $\H_{\sig,\tau,\eps}^m$-norm over any {\em finite-time} interval
of the form $[t_0,T]$ with $t_0>0$ (cf.\ Theorem \ref{ThmGP1}). More precisely, as a by-product of
Theorem \ref{ThmGP1}, we will prove that
\begin{equation*}
\lim_{\substack{\sig\to 0^+ \\ \tau\to 0^+   \\ \eps\to 0^+}}\sup_{t\in[t_0,T]}
\|\sed (t)z-\LL_{\sigma,\tau,\eps}\soo (t)\PP z\|_{\H_{\sig,\tau,\eps}^m}=0,
\end{equation*}
for every $R\geq 0$, $T>t_0>0$ and $z\in B_{\H_{\sig,\tau,\eps}^{2m+4}}(R)$.

The first three components of the solution
$\PP\sed (t)z$ are shown to converge to $\soo(t)\PP z$ in the
$\H^m_{0,0,0}$-norm on $[0,T]$, whereas
the history components $\eta^t$ and $\xi^t$ vanish on $[t_0,\infty]$
in the $\M_{\tau,\eps}^m$-norm and $\Q_\sig^{m+1}$-norm respectively,
due to the presence of possibly nonzero initial histories $\eta_0$ and $\xi_0$
(cf.\ Lemma \ref{etaestimm}).
\vskip2pt
\noindent
Besides, in the case $\tau=0$, the singular limit estimate
strengthens. Indeed, it turns out to hold on
{\em infinite-time} intervals far away from zero, uniformly with respect
to initial data lying inside any ball of $\H_{\sig,0,\eps}^{2m+4}$, namely
we get
\begin{equation*}
\lim_{\substack{\sig\to 0^+ \\ \eps\to 0^+}}
\sup_{z\in B_{\H_{\sig,0,\eps}^{2m+4}}(R)}\sup_{t\geq t_0}
\|S_{\sigma,0,\eps}(t)z-\LL_{\sigma,0,\eps}\soo (t)\PP z\|_{\H_{\sig,0,\eps}^m}=0,
\end{equation*}
for every $R\geq 0$ and $t_0>0$ (cf.\ Theorem \ref{ThmGP2}). Of course, to achieve
these results, the role played by the exponential stability of the limiting
semigroup $\soo(t)$ will be important.
\vskip2pt
\subsection{Some preliminary facts}
Before stating the main results of the section, we need a few preliminary results.
\begin{lemma}
\label{m-reg}
Let $m\geq 0$, $R\geq 0$ and $z\in B_{\H_{\sig,\tau,\eps}^m}(R)$. Then there exists $K_R\geq 0$
such that $\|\sed (t)z\|_{\H_{\sig,\tau,\eps}^m}\leq K_R$ for all $t\geq 0$.
\end{lemma}
\begin{proof}
By \eqref{etactrl} and \eqref{etactrl1}, it suffices to multiply
the equations of \ref{eq-mobile} by $u_t$ in $H^m$,
by $\vartheta$ in $H^m$, by $\eta$ in $\M^m_{\tau,\eps}$ and by $\xi$ in $\Q^{m+1}_\sig$
respectively and add the resulting equations.
\end{proof}

\noindent
The vanishing of the histories components $\eta^t$ and $\xi^t$
of $\sed(t)$ is issued in the following
\begin{lemma}
\label{etaestimm}
Let $m\geq 0$, $R\geq 0$ and $z=(u_0,u_1,\vartheta_0,\eta_0,\xi_0)\in B_{\H_{\sig,\tau,\eps}^m}(R)$.
Then there exists $K_R\geq 0$ such that the following facts hold:
\vskip5pt
\noindent
{\rm (a)} for every $\eps>0$ and $t\geq 0$,
\begin{equation}
\label{pezzo1}
\|\eta^t\|_{L^2_{\mu_\eps}(\R^+,H^{m+1})}\leq
\|\eta_0\|_{L^2_{\mu_\eps}(\R^+,H^{m+1})}e^{-\frac{\delta_1 t}{4\eps}}+K_R\sqrt{\eps};
\end{equation}
\vskip5pt
\noindent
{\rm (b)} for every $\tau>0$ and $t\geq 0$,
\begin{equation}
\label{tauhist}
\|\eta^t\|_{L^2_{\nu_\tau}(\R^+,H^m)}\leq
\|\eta_0\|_{L^2_{\nu_\tau}(\R^+,H^m)}e^{-\frac{\delta_3 t}{4}}+K_R\sqrt{\psi(\tau)};
\end{equation}
\vskip4pt
\noindent
{\rm (c)} for every $\sig>0$ and $t\geq 0$,
\begin{equation*}
\|\xi^t\|_{\Q_{\sig}^{m+1}}\leq
\|\xi_0\|_{\Q_{\sig}^{m+1}}e^{-\frac{\delta_2 t}{4\sig}}+K_R\sqrt{\sig}.
\end{equation*}
\end{lemma}
\begin{proof}
By arguing as in \cite[Lemma 5.4]{cps1} we immediately get assertions (a) and (c).
Let $C$ denote a generic positive constant depending on $R$
which may even vary from line to line within the same equation.
By multiplying the equation of $\eta$ by $\eta$ in $L^2_{\nu_\tau}(\R^+,H^m)$,
and taking \eqref{nutauips}, \eqref{nuip} and Lemma \ref{m-reg}
into account, we have
\begin{align*}
\frac{d}{dt}\|\eta\|^2_{L^2_{\nu_\tau}(\R^+,H^m)}+\delta_3\|\eta\|^2_{L^2_{\nu_\tau}(\R^+,H^m)}
&\leq C\displaystyle\int_0^\infty \nu_\tau(s)\|A^{m/2}\eta(s)\|ds\\
\noalign{\vskip2pt}
&\leq C\Big(\int_0^\infty\nu_\tau(s)ds\Big)^{1/2}
\Big(\int_0^\infty \nu_\tau(s)\|A^{m/2}\eta(s)\|^2ds\Big)^{1/2} \\
\noalign{\vskip2pt}
&\leq C\sqrt{\psi(\tau)}\|\eta\|_{L^2_{\nu_\tau}(\R^+,H^m)}
\leq \frac{\delta_3}{2}\|\eta\|_{L^2_{\nu_\tau}(\R^+,H^m)}^2
+C\psi(\tau),
\end{align*}
so that, by the Gronwall Lemma, we immediately obtain (b).
\end{proof}

\begin{definition}
For every $m\geq 0$, $\eta_0\in \M_{\tau,\eps}^m$ and $\xi_0\in\Q_{\sig}^{m+1}$,
let us set for every $t\geq 0$
\begin{equation}
\label{defups}
\Upsilon^m_{\sig,\tau,\eps}(t)=
\|\eta_0\|_{L^2_{\mu_\eps}(\R^+,H^{m+1})}e^{-\frac{\delta_1 t}{4\eps}}+
\|\eta_0\|_{L^2_{\nu_\tau}(\R^+,H^m)}e^{-\frac{\delta_3 t}{4}}+
\|\xi_0\|_{\Q_{\sig}^{m+1}}e^{-\frac{\delta_2 t}{4\sig}}.
\end{equation}
Furthermore, we introduce the maps $\Pi_\flat:[0,1]^3\to\R^+$ and $\Pi_\sharp:[0,1]\to\R^+$,
\begin{align*}
\Pi_\flat(\sig,\tau,\eps)&=\sqrt[4]{\eps}+\sqrt[4]{\sigma}+\sqrt[4]{\psi(\tau)}, \\
\noalign{\vskip4pt}
\Pi_\sharp(\tau)&=\sqrt{\psi(\tau)}+\sqrt{\phi(\tau)}.
\end{align*}
Observe that $\Pi_\flat$ and $\Pi_\sharp$ are continuous
with $\Pi_\flat(0,0,0)=\Pi_\sharp(0)=0$.
\end{definition}

\begin{proposition}
\label{azeero}
For every $m\geq 0$, $R\geq 0$, $\eta_0\in B_{\M_{\tau,\eps}^m}(R)$,
$\xi_0\in B_{\Q_{\sig}^{m+1}}(R)$ and $t_0>0$
\begin{equation*}
\lim_{\substack{\sig\to 0^+ \\  \tau\to 0^+  \\ \eps\to 0^+}}\sup_{t\geq t_0}
\Upsilon^m_{\sig,\tau,\eps}(t)=0.
\end{equation*}
\end{proposition}
\begin{proof}
The first and third summands of $\Upsilon^m_{\sig,\tau,\eps}$ vanish exponentially. Moreover,
observe that, by \eqref{nuip} and since $\{\tau\mapsto\nu_\tau\}$ is increasing,
$\|\eta_0\|_{L^2_{\nu_\tau}(\R^+,H^m)}$ converges to zero by the Monotone
Convergence Theorem.
\end{proof}

\subsection{Case $\boldsymbol{\tau>0}$ : the convergence estimate}
\label{convres}
We are now ready to state the main result of the section, which gives
a convergence estimate of $\sed(t)$ towards $\soo(t)$ in the
norm of $\H_{\sig,\tau,\eps}^m$, for any $m\geq 0$, over finite-time intervals.

\begin{theorem}
\label{ThmGP1}
For every $m\geq 0$, $R\geq 0$, $T>0$ and $z\in B_{\H_{\sig,\tau,\eps}^{2m+4}}(R)$,
there exist two constants $K_{R}\geq 0$ and $Q_{R,T}\geq 0$ such that
\begin{equation*}
\|\sed (t)z-\LL_{\sigma,\tau,\eps}\soo (t)\PP z\|_{\H_{\sig,\tau,\eps}^m}
\leq \Upsilon^m_{\sig,\tau,\eps}(t)+K_R\Pi_\flat(\sig,\tau,\eps)+Q_{R,T}\Pi_\sharp(\tau),
\end{equation*}
for every $t\in [0,T]$.
\end{theorem}

\begin{remark}
The regularity assumption on the initial data $z$ could be relaxed to get
a rougher convergence estimate on finite-time intervals. On the
other hand, the price one has to pay is that also the constant $K_{R}$ which appears
in the above theorem would depend on the time interval. With the higher regularity that
we require, instead, we are able to exploit the exponential stability of the
limiting semigroup $\soo(t)$ and to have, at least in the case $\tau=0$, the
convergence estimate holding uniformly in time. So, for $m=0$, we get
a convergence estimate for the thermoviscoelastic model starting
with initial data having {\em four} levels of regularity above the regularity of the base phase-space.
In the thermoelastic plate model considered in \cite{GS1} (essentially, w.r.t.\ our notation,
the case when $\sigma=\tau=0$) one needs to go just {\em two} levels of regularity above.
Finally, in the case of single parabolic and hyperbolic equations with memory (cf.\ \cite{cps1,cps2})
it suffices to require {\em one} level of regularity above the basic regularity
to get a time dependent control.
\end{remark}

\begin{proof}[Proof of Theorem \ref{ThmGP1}]
Let $m\geq 0$, $R\geq 0$ and $z=(u_0,u_1,\vartheta_0,\eta_0,\xi_0)
\in B_{\H_{\sig,\tau,\eps}^{2m+4}}(R)$. Since
$$
\sed(t)z=(\PP\sed(t)z,\QQ_{\tau,\eps}\sed(t)z,\QQ_\sigma\sed(t)z),\qquad t\geq 0,
$$
we get the assertion if we prove that, for $T>0$, there exist
$K_{R}\geq 0$ and $Q_{R,T}\geq 0$ with
\begin{align}
\label{un}
\|\PP\sed (t)z-\soo (t)\PP z\|_{\H_{0,0,0}^m}
& \leq K_R\,\Pi_\flat(\sig,\tau,\eps)+Q_{R,T}\Pi_\sharp(\tau)  \\
\noalign{\vskip3pt}
\label{du}
\|\QQ_{\tau,\eps}\sed (t)z\|_{\M_{\tau,\eps}^m}+\|\QQ_{\sigma}\sed (t)z\|_{\Q_\sigma^{m+1}}&
\leq\Upsilon^m_{\sigma,\tau,\eps}(t)+K_R\big(\sqrt{\eps}+\sqrt{\sig}+\sqrt{\psi(\tau)}\big),
\end{align}
for every $t\in[0,T]$. By combining inequalities (a), (b) and (c)
of Lemma \ref{etaestimm}, we immediately obtain \eqref{du}.
Then, we turn to the proof of inequality \eqref{un}. Let us set
\begin{align*}
\bar u(t)&= \hat u(t)-u(t), \\
\noalign{\vskip3pt}
\bar u_t(t)&= \hat u_t(t)-u_t(t), \\
\noalign{\vskip3pt}
\bar \vartheta(t) &= \hat \vartheta(t)-\vartheta(t), \\
\noalign{\vskip3pt}
\bar\eta^t&=\hat\eta^t-\eta^t, \\
\noalign{\vskip3pt}
\bar\xi^t&=\hat\xi^t-\xi^t,
\end{align*}
where $(\hat u,\hat u_t,\hat \vartheta,\hat\eta,\hat\xi)$
denotes the solution to \ref{eq-mobile} with initial data $z$,
while $(u,u_t,\vartheta)$ stands for the solution to \ref{eq-limite}
with initial data $\PP z$. Besides, $\eta^t$ (resp.\ $\xi^t$) denotes the solution at
time $t$ of the Cauchy problem in $\M_{\tau,\eps}^0$ (resp.\ $\Q_\sigma^1$)
$$
\begin{cases}
\eta_t=T_{\tau,\eps}\eta+\vartheta,\qquad t>0,\\
\noalign{\vskip2pt}
\eta^0=\eta_0,
\end{cases}
\,\,\qquad\,\,
\begin{cases}
\xi_t=T_\sigma\xi+u_t,\qquad t>0,\\
\noalign{\vskip2pt}
\xi^0=\xi_0.
\end{cases}
$$
These problems reconstruct the missing components of the limiting semigroup $\soo(t)$
which are needed in order to perform the comparison argument (cf.\ \cite{cps1,cps2}).
Then, it can be readily checked that $(\bar u,\bar u_t,\bar\vartheta,\bar\eta,\bar\xi)$
solves the system
$$
\begin{cases}
\displaystyle
\bar u_{tt}+A^2\bar u+\int_0^\infty \beta_\sigma(s)A^2\hat\xi(s)ds-A^2u_t-A\bar\vartheta=0, \\
\noalign{\vskip4pt}
\displaystyle
\bar\vartheta_t+\phi(\tau)\bar\vartheta+\phi(\tau)\vartheta
+\int_0^\infty \nu_\tau(s)\hat\eta(s)ds
+\displaystyle\int_0^\infty \mu_\eps(s)A\hat\eta(s)ds-A\vartheta+A\bar u_t=0, \\
\noalign{\vskip5pt}
\bar\eta_t=T_{\tau,\eps}\bar\eta+\bar\vartheta, \\
\noalign{\vskip5pt}
\bar\xi_t=T_\sigma\bar\xi+\bar u_t, \\
\noalign{\vskip5pt}
(\bar u(0),\bar u_t(0),\bar\vartheta(0),\bar\eta^0,\bar\xi^0)=(0,0,0,0,0).
\end{cases}
$$
By multiplying the first equation by $\bar u_t$ in $H^m$, the second
by $\bar\vartheta$ in $H^m$, the third by $\bar\eta$ in $\M_{\tau,\eps}^m$
and the fourth by $\bar\xi$ in $\Q_\sig^{m+1}$ we obtain, respectively,
\begin{align*}
{\textstyle\frac{1}{2}}\frac{d}{dt}\big(\|A^{(m+2)/2}\bar u\|^2+\|A^{m/2}\bar u_t\|^2\big)
+\l \hat\xi,\bar u_t\r_{\Q_\sig^{m+1}}  \\
\noalign{\vskip2pt}
-\l A^{(m+2)/2}u_t,A^{(m+2)/2}\bar u_t\r
-\langle A^{(m+1)/2}\bar \vartheta,A^{(m+1)/2}\bar u_t\rangle &=0, \\
\noalign{\vskip7pt}
{\textstyle\frac{1}{2}}\frac{d}{dt}\|A^{m/2}\bar \vartheta\|^2
+\phi(\tau)\|A^{m/2}\bar\vartheta\|^2
+\phi(\tau)\l A^{m/2}\vartheta, A^{m/2}\bar\vartheta\r \\
\noalign{\vskip2pt}
+\langle\hat\eta,\bar\vartheta \rangle_{\M_{\tau,\eps}^m}
-\langle A^{(m+1)/2}\vartheta,
A^{(m+1)/2}\bar\vartheta\rangle+\langle A^{(m+1)/2}\bar u_t,A^{(m+1)/2}\bar\vartheta\rangle &=0, \\
\noalign{\vskip7pt}
{\textstyle\frac{1}{2}}\frac{d}{dt}\|\bar\eta\|^2_{\M_{\tau,\eps}^m}-
\langle T_{\tau,\eps}\bar\eta,\bar\eta\rangle_{\M_{\tau,\eps}^m}-
\langle\bar\eta,\bar\vartheta \rangle_{\M_{\tau,\eps}^m} &=0, \\
\noalign{\vskip6pt}
{\textstyle\frac{1}{2}}\frac{d}{dt}\|\bar\xi\|^2_{\Q_\sig^{m+1}}-
\langle T_\sig\bar\xi,\bar\xi\rangle_{\Q_\sig^{m+1}}-
\langle\bar\xi,\bar u_t \rangle_{\Q_\sig^{m+1}} &=0.
\end{align*}
Taking \eqref{etactrl}-\eqref{etactrl1} into account, and adding the above
identities, we end up with
\begin{equation*}
\frac{d}{dt}\big(\|A^{(m+2)/2}\bar u\|^2+\|A^{m/2}\bar u_t\|^2+\|A^{m/2}\bar \vartheta\|^2
+\|\bar\eta\|^2_{\M_{\tau,\eps}^m}+\|\bar\xi\|^2_{\Q_\sig^{m+1}}\big)\leq 2I_\eps+2J_\sigma+2K_\tau,
\end{equation*}
where we have set
\begin{align*}
I_\eps(t)&=-\int_0^\infty \mu_\eps(s)\langle A^{(m+1)/2}\eta^t(s),A^{(m+1)/2}\bar\vartheta(t)\rangle ds
+\langle A^{(m+1)/2}\vartheta(t),A^{(m+1)/2}\bar\vartheta(t)\rangle, \\
J_\sig(t)&=-\int_0^\infty \beta_\sigma(s)\langle A^{(m+2)/2}\xi^t(s),A^{(m+2)/2}\bar u_t(t)\rangle ds
+\langle A^{(m+2)/2}u_t(t),A^{(m+2)/2}\bar u_t(t)\rangle, \\
\noalign{\vskip4pt}
K_\tau(t)&=-\int_0^\infty \nu_\tau(s)\langle A^{m/2}\eta^t(s),A^{m/2}\bar\vartheta(t)\rangle ds
-\phi(\tau)\langle A^{m/2}\vartheta(t),A^{m/2}\bar\vartheta(t)\rangle.
\end{align*}
We stress, for later use, that the above term $K_\tau(t)$ appears under
the assumption that $\tau>0$, whereas we would simply have
$K_0(t)=0$ for all $t\geq 0$ in the case $\tau=0$,
since $\nu_0=0$ and $\phi(0)=0$. We shall write,
for all $t\geq 0$, $I_\eps(t)=\sum_{j=1}^5 I_j(t)$ and
$J_\sig(t)=\sum_{j=1}^5 J_j(t)$, being the $I_j$s and the
$J_j$s defined, respectively, by
\begin{align*}
I_1(t)&=\int_{\sqrt{\eps}}^\infty s\mu_\eps(s)
\bra{A^{(m+1)/2}\vartheta(t),A^{(m+1)/2}\bar \vartheta(t)}ds,\\
\noalign{\vskip4pt}
I_2(t)&=-\int_{\sqrt{\eps}}^\infty
       \mu_\eps(s)\bra{A^{(m+1)/2}\eta^t(s),A^{(m+1)/2}\bar \vartheta(t)}ds,\\
\noalign{\vskip4pt}
I_3(t)&=-\int_{\min\{\sqrt{\eps},t\}}^{\sqrt{\eps}}\mu_\eps(s)
       \bra{A^{(m+1)/2}\eta_0(s-t),A^{(m+1)/2}\bar \vartheta(t)} ds,\\
\noalign{\vskip4pt}
I_4(t)&=\int_{\min\{\sqrt{\eps},t\}}^{\sqrt{\eps}}
      (s-t)\mu_\eps(s)\bra{A^{(m+1)/2}\vartheta(t),A^{(m+1)/2}\bar\vartheta(t)} ds,\\
\noalign{\vskip4pt}
I_5(t)&=\int_0^{\sqrt{\eps}} \mu_\eps(s)\bigg[
      \int_0^{\min\{s,t\}}\bra{A^{(m+1)/2}(\vartheta(t)-\vartheta(t-y)),
      A^{(m+1)/2}\bar \vartheta(t)}dy\bigg]ds,
\end{align*}
and
\begin{align*}
J_1(t)&=\int_{\sqrt{\sig}}^\infty s\beta_\sigma(s)
\bra{A^{(m+2)/2}u_t(t),A^{(m+2)/2}\bar u_t(t)}ds,\\
\noalign{\vskip4pt}
J_2(t)&=-\int_{\sqrt{\sig}}^\infty
       \beta_\sigma(s)\bra{A^{(m+2)/2}\xi^t(s),A^{(m+2)/2}\bar u_t(t)}ds,\\
\noalign{\vskip4pt}
J_3(t)&=-\int_{\min\{\sqrt{\sig},t\}}^{\sqrt{\sig}}\beta_\sigma(s)
       \bra{A^{(m+2)/2}\xi_0(s-t),A^{(m+2)/2}\bar u_t(t)} ds,\\
\noalign{\vskip4pt}
J_4(t)&=\int_{\min\{\sqrt{\sig},t\}}^{\sqrt{\sig}}
      (s-t)\beta_\sigma(s)\bra{A^{(m+2)/2}u_t(t), A^{(m+2)/2}\bar u_t(t)} ds,\\
\noalign{\vskip4pt}
J_5(t)&=\int_0^{\sqrt{\sig}} \beta_\sigma(s)\bigg[
      \int_0^{\min\{s,t\}}\bra{A^{(m+2)/2}(u_t(t)-u_t(t-y)),
      A^{(m+2)/2}\bar u_t(t)}dy\bigg]ds.
\end{align*}
In the following, we shall denote by $C\geq 0$ a generic constant
which may even vary from line to line and may depend on $R$, but
it is independent of $\sig$, $\tau$ and $\eps$.
By virtue of Lemma \ref{m-reg}, we have $\|\sed (t)z\|_{\H_{\sig,
\tau,\eps}^{2m+4}}\leq C$ for all $t\geq 0$. In particular,
\begin{align}
\label{2-control-a}
\|A^{m+3}\hat u(t)\|+\|A^{m+2}\hat u_t(t)\|+\|A^{m+2}\hat \vartheta(t)\|+
\|\hat\eta^t\|_{\M^{2m+4}_{\tau,\eps}}+\|\hat\xi^t\|_{\Q^{2m+5}_\sigma}  &\leq C, \\
\noalign{\vskip3pt}
\label{2-control-b}
\|A^{m+3}u(t)\|+\|A^{m+2}u_t(t)\|+\|A^{m+2}\vartheta(t)\|+
\|\eta^t\|_{\M^{2m+4}_{\tau,\eps}}+\|\xi^t\|_{\Q^{2m+5}_\sigma}&\leq C,
\end{align}
for all $t\geq 0$. Furthermore, since $\soo(t)$ is exponentially stable,
there exists $\varpi>0$ with
\begin{equation}
\label{limit-stabil}
\|Au(t)\|+\|u_t(t)\|+\|\vartheta(t)\| \leq Ce^{-\varpi t},\qquad\forall t\geq 0.
\end{equation}
Concerning the treatment of the terms $I_j$s and $J_j$s, we will proceed
on the line of \cite{cps1,cps2} but strengthening the estimates, whenever it is possible, through
the first order energy exponential decay furnished by \eqref{limit-stabil}. Observe
first that, due to \eqref{K3}-\eqref{K32} and \eqref{riscke},
\begin{equation}
\label{intmueps1}
\int_{\sqrt{\eps}}^\infty s\mu_\eps(s)ds \leq C\eps,\quad\forall\,\eps>0, \qquad
\int_{\sqrt{\sig}}^\infty s\beta_\sigma(s)ds \leq C\sig,\quad\forall\,\sig>0.
\end{equation}
Hence, by \eqref{2-control-a}, \eqref{2-control-b}, \eqref{limit-stabil}
and \eqref{intmueps1}, we immediately get
\begin{align*}
I_1(t)&=\int_{\sqrt{\eps}}^\infty s\mu_\eps(s)
\bra{\vartheta(t),A^{m+1}\bar \vartheta(t)}ds  \\
\noalign{\vskip2pt}
&\leq C\eps\|A^{m+1}\bar\vartheta(t)\|\|\vartheta(t)\|\leq C\eps e^{-\varpi t},
\qquad\forall\, t\geq 0,\\
\noalign{\vskip5pt}
J_1(t)&=\int_{\sqrt{\sig}}^\infty s\beta_\sigma(s)
\bra{u_t(t),A^{m+2}\bar u_t(t)}ds \\
\noalign{\vskip2pt}
&\leq C\sig\|A^{m+2}\bar u_t(t)\|\|u_t(t)\|\leq C\sigma e^{-\varpi t},
\qquad\forall\, t\geq 0.
\end{align*}
Let us now prove that there holds
\begin{align}
\label{lim-sto}
\|\eta^t\|_{L^2_{\mu_\eps}(\R^+,H^{m+1})}^2 &\leq
Ce^{-\frac{\delta_1 t}{\eps}}
+C\sqrt{\eps}\,e^{-\varpi t},
\qquad\forall\, t\geq 0, \\
\noalign{\vskip1pt}
\label{lim-sto2}
\|\xi^t\|_{\Q_\sig^{m+1}}^2 &\leq
Ce^{-\frac{\delta_2 t}{\sig}}+C\sqrt{\sig}\,e^{-\varpi t},
\qquad\forall\, t\geq 0,
\end{align}
Indeed, arguing as in \cite[Lemma 5.4]{cps1}, we readily obtain
\begin{align*}
\|\eta^t\|_{L^2_{\mu_\eps}(\R^+,H^{2m+2})} &\leq Ce^{-\frac{\delta_1t}{4\eps}}+C\sqrt{\eps},
\qquad\forall\, t\geq 0, \\
\noalign{\vskip1pt}
\|\xi^t\|_{\Q_\sigma^{2m+3}} &\leq Ce^{-\frac{\delta_2 t}{4\sig}}+C\sqrt{\sig},
\qquad\forall\, t\geq 0.
\end{align*}
Whence, by multiplying the equation for $\eta$
times $\eta$ in $L^2_{\mu_\eps}(\R^+,H^{m+1})$, in light of \eqref{limit-stabil},
\begin{align*}
\frac{d}{dt}\|\eta^t\|^2_{L^2_{\mu_\eps}(\R^+,H^{m+1})}
&+\textstyle{\frac{\delta_1}{\eps}}\|\eta^t\|^2_{L^2_{\mu_\eps}(\R^+,H^{m+1})}
\leq 2\|\vartheta(t)\|\displaystyle\int_0^\infty \mu_\eps(s)\|A^{m+1}\eta^t(s)\|ds \\
\noalign{\vskip1pt}
&\leq \textstyle{\frac{2}{\sqrt{\eps}}}\|\vartheta(t)\|\|\eta^t\|^2_{L^2_{\mu_\eps}(\R^+,H^{2m+2})}
\leq  \textstyle{\frac{C}{\sqrt{\eps}}}e^{-(\varpi+\frac{\delta_1}{4\eps}) t}
+Ce^{-\varpi t},
\end{align*}
which yields \eqref{lim-sto} via the Gronwall Lemma.
In a similar fashion, again by \eqref{limit-stabil}, we get
\begin{align*}
\frac{d}{dt}\|\xi^t\|^2_{\Q_\sig^{m+1}}
+\textstyle{\frac{\delta_2}{\sig}} \|\xi^t\|^2_{\Q_\sig^{m+1}}
&\leq 2\|u_t(t)\|\int_0^\infty \beta_\sigma(s)\|A^{m+2}\xi^t(s)\|ds\\
&\leq \textstyle{\frac{2}{\sqrt{\sig}}}\|u_t(t)\|\|\xi^t\|_{\Q_\sig^{2m+3}}
\leq  \textstyle{\frac{C}{\sqrt{\sig}}} e^{-(\varpi+\frac{\delta_2}{4\sig})t}+Ce^{-\varpi t},
\end{align*}
which yields inequality \eqref{lim-sto2}.
By means of \eqref{2-control-a}-\eqref{2-control-b},
\eqref{intmueps1} and \eqref{lim-sto}-\eqref{lim-sto2} we have
\begin{align*}
I_2(t) &\leq C\int_{\sqrt{\eps}}^\infty \mu_\eps(s)\|A^{(m+1)/2}\eta^t(s)\|ds \\
&\leq C\sqrt{\eps}\|\eta^t\|_{L^2_{\mu_\eps}(\R^+,H^{m+1})}
\leq C\sqrt{\eps}\,e^{-\frac{\delta_1 t}{2\eps}}
+C\sqrt{\eps}\, e^{-\frac{\varpi}{2} t},\qquad\forall\, t\geq 0, \\
\noalign{\vskip8pt}
J_2(t) &\leq C\int_{\sqrt{\sig}}^\infty \beta_\sigma(s)\|A^{(m+2)/2}\xi^t(s)\|ds \\
&\leq C\sqrt{\sig}\|\xi^t\|_{\Q_\sig^{m+1}}
\leq C\sqrt{\sig}\,e^{-\frac{\delta_2 t}{2\sig}}+C\sqrt{\sig}\, e^{-\frac{\varpi}{2} t},
\qquad\forall\, t\geq 0.
\end{align*}
Taking \eqref{K3}-\eqref{K32}, \eqref{normalize1}-\eqref{normalize3} and
\eqref{2-control-a}-\eqref{2-control-b} into account, we get, for $t<\sqrt{\eps}$,
\begin{align*}
I_3(t)
&\leq C\int_{t}^{\sqrt{\eps}}\mu_\eps(s)
\|A^{(m+1)/2}\eta_0(s-t)\| ds\\
&\leq
Ce^{-\frac{\delta_1 t}{\eps}}\Big(\int_{0}^{\infty}\mu_\eps(s)
ds\Big)^{1/2}\|\eta_0\|_{L^2_{\mu_\eps}(\R^+,H^{m+1})}\leq\textstyle{\frac{C}{\sqrt{\eps}}}\,
e^{-\frac{\delta_1 t}{\eps}}, \qquad\forall\, t\geq 0,\\
\noalign{\vskip3pt}
J_3(t)
&\leq C\int_{t}^{\sqrt{\sig}}\beta_\sigma(s)\|A^{(m+2)/2}\xi_0(s-t)\| ds\\
&\leq
Ce^{-\frac{\delta_2 t}{\sig}}\Big(\int_{0}^{\infty}\beta_\sigma(s)
ds\Big)^{1/2}\|\xi_0\|_{\Q_\sig^{m+1}}\leq\textstyle{\frac{C}{\sqrt{\sig}}}\,
e^{-\frac{\delta_2 t}{\sig}},\qquad\forall\, t\geq 0.
\end{align*}
Arguing in a similar fashion, there holds
\begin{align*}
I_4(t) &\leq C e^{-\frac{\delta_1 t}{\eps}}
\int_{0}^{\infty}s\mu_\eps(s)ds
= C e^{-\frac{\delta_1 t}{\eps}},  \qquad\forall\, t\geq 0,\\
\noalign{\vskip3pt}
J_4(t) &\leq C e^{-\frac{\delta_2 t}{\sig}}
\int_{0}^{\infty}s\beta_\sigma(s)ds
= C e^{-\frac{\delta_2 t}{\sig}},\qquad\forall\, t\geq 0.
\end{align*}
Observe now that
\begin{equation}
\label{controlDer}
\|\partial_t S_{0,0,0}(t)\PP z\|_{\H_{0,0,0}^0}\leq C,\qquad\forall\, t\geq 0.
\end{equation}
Indeed, from the equations of \ref{eq-limite}, by virtue
of \eqref{2-control-b},
\begin{align*}
\|\vartheta_t(t)\|& \leq \|A\vartheta(t)\|+\|Au_t(t)\|\leq C, \qquad\forall\, t\geq 0,\\
\noalign{\vskip4pt}
\|u_{tt}(t)\| &\leq\|A^2u(t)\|+\|A^2u_t(t)\|+\|A\vartheta(t)\|\leq C,\qquad\forall\, t\geq 0,
\end{align*}
so that \eqref{controlDer} readily follows. In particular \eqref{controlDer} yields
\begin{align*}
\label{mod-cont-th}
\|u_t(t)-u_t(t-y)\|+\|\vartheta(t)-\vartheta(t-y)\| &\leq \|\soo(t-y)(\soo(y)
\PP z-\PP z)\|_{\H_{0,0,0}^0} \\
& \leq Ce^{-\varpi t}\int_{0}^y\|\partial_t S_{0,0,0}(\varsigma)\PP z
\|_{\H_{0,0,0}^0}d\varsigma \\
&\leq Ce^{-\varpi t}y,
\notag
\end{align*}
for every $t\geq 0$ and $y\in[0,t]$. Hence, by
\eqref{normalize1}-\eqref{normalize3} and
\eqref{2-control-a}-\eqref{2-control-b}, we obtain
\begin{align*}
I_5(t) & \leq \|A^{m+1}\bar \vartheta(t)\|\int_0^{\sqrt{\eps}} \mu_\eps(s)
\int_0^{\min\{s,t\}}\|\vartheta(t)-\vartheta(t-y)\|dyds \\
\noalign{\vskip1pt}
& \leq C\sqrt{\eps}\,e^{-\varpi t},\qquad\forall\, t\geq 0, \\
\noalign{\vskip4pt}
J_5(t) & \leq \|A^{m+2}\bar u_t(t)\|\int_0^{\sqrt{\sig}} \beta_\sigma(s)
\int_0^{\min\{s,t\}}\|u_t(t)-u_t(t-y)\|dyds \\
\noalign{\vskip1pt}
&\leq C\sqrt{\sig}\,e^{-\varpi t},\qquad\forall\, t\geq 0.
\end{align*}
We now turn to the estimate of $K_\tau$. Taking
condition \eqref{nuip} as well as (b) of Lemma \ref{etaestimm}
into account, we obtain
\begin{align*}
K_\tau(t)&=-\int_0^\infty \nu_\tau(s)\langle A^{m/2}\eta^t(s),A^{m/2}\bar\vartheta(t)\rangle ds
-\phi(\tau)\langle A^{m/2}\vartheta(t),A^{m/2}\bar\vartheta(t)\rangle,  \\
\noalign{\vskip2pt}
& \leq C\displaystyle\int_0^\infty \nu_\tau(s)\|A^{m/2}\eta^t(s)\|ds+C \phi(\tau) \\
\noalign{\vskip4pt}
&\leq C\Big(\int_0^\infty\nu_\tau(s)ds\Big)^{1/2}
\Big(\int_0^\infty \nu_\tau(s)\|A^{m/2}\eta(s)\|^2ds\Big)^{1/2}+C\phi(\tau) \\
\noalign{\vskip6pt}
&\leq C\sqrt{\psi(\tau)}\|\eta_0\|_{L^2_{\nu_\tau}(\R^+,H^m)}
e^{-\frac{\delta_3 t}{4}}+C\psi(\tau)+C\phi(\tau) \\
\noalign{\vskip10pt}
&\leq C\sqrt{\psi(\tau)}e^{-\frac{\delta_3 t}{4}}+C\psi(\tau)+C\phi(\tau),
\qquad\forall\, t\geq 0.
\end{align*}
Therefore, collecting the previous inequalities, we end up with
\begin{equation*}
\frac{d}{dt}\|\PP\sed (t)z-\soo (t)\PP z\|_{\H_{0,0,0}^m}^2
\leq \varphi_1(t)+\varphi_2(t),\qquad\forall\, t\geq 0,
\end{equation*}
where we have set
\begin{align*}
\varphi_1(t)&=C\big[(\sqrt{\eps}+\sqrt{\sig})e^{-\frac{\varpi}{2} t}
+\textstyle{\frac{1}{\sqrt{\eps}}}\,e^{-\frac{\delta_1 t}{2\eps}}
+\textstyle{\frac{1}{\sqrt{\sig}}}\,e^{-\frac{\delta_2 t}{2\sig}}\big], \\
\noalign{\vskip5pt}
\varphi_2(t)&=C\sqrt{\psi(\tau)}e^{-\frac{\delta_3 t}{4}}+C\psi(\tau)+C\phi(\tau).
\end{align*}
Notice that, there holds
\begin{align*}
\int_0^t \varphi_1(\varsigma)d\varsigma &\leq C(\sqrt{\eps}+\sqrt{\sig}),\qquad\forall t\geq 0, \\
\noalign{\vskip2pt}
\int_0^t \varphi_2(\varsigma)d\varsigma &\leq C\sqrt{\psi(\tau)}
+C(\psi(\tau)+\phi(\tau))t,
\qquad\forall\, t\geq 0.
\end{align*}
Consequently, by integrating the above differential inequality in time,
we can find two constants $K_R\geq 0$ and $Q_{R,T}\geq 0$ such that
\begin{equation*}
\|(\bar u(t),\bar u_t(t),\bar\vartheta(t))\|_{\H_{0,0,0}^m}
\leq K_R\,\Pi_\flat(\sig,\tau,\eps)+Q_{R,T}\Pi_\sharp(\tau),\qquad\forall\, t\geq 0,
\end{equation*}
which proves \eqref{un}. The proof is now complete.
\end{proof}

\subsection{Case $\boldsymbol{\tau=0}$ : uniform in time estimate}
As a straightforward but important corollary of the main
Theorem \ref{ThmGP1}, in the case $\tau=0$, we obtain the following

\begin{theorem}
\label{ThmGP2}
For every $m\geq 0$, $R\geq 0$ and $z\in B_{\H_{\sig,0,\eps}^{2m+4}}(R)$
there exists $K_{R}\geq 0$ with
\begin{equation*}
\|S_{\sigma,0,\eps}(t)z-\LL_{\sigma,0,\eps}\soo (t)\PP z\|_{\H_{\sig,0,\eps}^m}
\leq \|\eta_0\|_{\M_{0,\eps}^m}e^{-\frac{\delta_1 t}{4\eps}}+
\|\xi_0\|_{\Q_{\sig}^{m+1}}e^{-\frac{\delta_2 t}{4\sig}}+
K_R(\sqrt[4]{\eps}+\sqrt[4]{\sigma}),
\end{equation*}
for every $t\in [0,\infty)$.
\end{theorem}
\begin{proof}
It suffices to retrace the steps in the proof of Theorem \ref{ThmGP1},
taking into account that $\Pi_\flat(\sig,0,\eps)=\sqrt[4]{\eps}+\sqrt[4]{\sigma}$
and $K_0(t)=0$ for every $t\geq 0$, which in turn yields $Q_{R,T}=0$.
\end{proof}

\bigskip

\section*{Appendix: failure of exponential decay}
\setcounter{equation}{0} \setcounter{subsection}{0}
\renewcommand{\theequation}{A.\arabic{equation}}
\renewcommand{\thesubsection}{A.\arabic{subsection}}

\theoremstyle{plain}
\newtheorem{theoremA}{Theorem}[section]
\renewcommand{\thetheoremA}{A.\arabic{theoremA}}
\newtheorem{lemmaA}[theoremA]{Lemma}
\renewcommand{\thelemmaA}{A.\arabic{lemmaA}}
\newtheorem{propositionA}[theoremA]{Proposition}
\renewcommand{\thepropositionA}{A.\arabic{propositionA}}
\newtheorem{corollaryA}[theoremA]{Corollary}
\renewcommand{\thecorollaryA}{A.\arabic{corollaryA}}
\theoremstyle{definition}
\newtheorem{definitionA}[theoremA]{Definition}
\renewcommand{\thedefinitionA}{A.\arabic{definitionA}}
\theoremstyle{remark}
\newtheorem{remarkA}[theoremA]{Remark}
\renewcommand{\theremarkA}{A.\arabic{remarkA}}

\noindent
Here we want to consider some variants of our model with no energy relaxation (i.e., $a\equiv 0$), in order
to show that the relaxation of the heat flux and/or of the strain
may not ensure the exponential stability of the corresponding semigroup.

Let $\Omega\subset\R^2$ be a smooth bounded domain and let $\alpha,\sigma\geq 0$.
In this section we consider the thermoelastic system with memory in an abstract setting
\begin{equation*}
\begin{cases}
u_{tt}+\displaystyle\int_0^\infty h'(s)A^2u(t-s)ds+A^2u-A^\sigma\vartheta=0, \\
\noalign{\vskip7pt}
\vartheta_t+\displaystyle\int_0^\infty k(s)A^\alpha\vartheta(s-t)ds+A^\sigma u_t=0.
\end{cases}
\end{equation*}
Notice that the corresponding memory free model, studied, e.g.,\ in \cite{riv-rak2},
\begin{equation*}
\begin{cases}
u_{tt}+A^2u-A^\sigma\vartheta=0, \\
\noalign{\vskip7pt}
\vartheta_t+A^\alpha\vartheta+A^\sigma u_t=0,
\end{cases}
\end{equation*}
includes, for $A=-\Delta$, as particular cases:
\vskip4pt
\noindent
- {\em thermoelastic plates}, for $\alpha=\sigma=1$;
\vskip2pt
\noindent
- {\em viscoelasticity}, for $\alpha=0$ and $\sig=1$.

\subsection{Preliminaries and main results}
Let $k:\R^+\to\R^+$ and $h:\R^+\to\R^+$ be smooth, decreasing,
summable functions and set $\mu(s)=-k'(s)$ and $\beta(s)=-h'(s)$, where
\begin{align}
\label{K0A}
&\mu,\beta \in C^1(\R^+)\cap L^1(\R^+), \\
\noalign{\vskip2pt}
\label{K1A}
& \mu(s)\ge 0, \quad\beta(s)\geq 0,
\qquad \forall\, s\in\R^+,\\
\noalign{\vskip4pt}
\label{K2A}
& \mu'(s)\le 0, \quad\beta'(s)\leq 0,
\qquad \forall\, s\in\R^+, \\
\noalign{\vskip2pt}
\label{add-int}
& 0<\int_0^\infty
s^2\mu(s)ds<\infty,\qquad \int_0^\infty \beta(s)ds>0.
\end{align}
Moreover, we assume that $A$ is a (strictly) positive selfadjoint linear operator on
$L^2(\Omega)$ with domain $\D(A)$ which admits a diverging sequence of positive
eigenvalues $\{\gamma_n\}_{n\geq 1}$.
\vskip2pt
We introduce the scale of Hilbert spaces
$H^\alpha=\D(A^{\alpha/2})$, $\alpha\in\R$, endowed with the inner
products $\langle u_1,u_2\rangle_ {H^\alpha}=\langle
A^{\alpha/2}u_1,A^{\alpha/2}u_2\rangle$ and we consider the
weighted spaces
$$
\M_{\alpha}=L^2_{\mu}(\R^+,H^{\alpha}),\qquad
\Q=L^2_{\beta}(\R^+,H^2),\qquad \alpha\in\R,
$$
endowed, respectively, with the inner products
\begin{equation*}
\l\eta_1,\eta_2\r_{\M_\alpha}
=\int_0^\infty\mu(s)\langle \eta_1(s),\eta_2(s)\rangle_ {H^\alpha}ds,
\qquad
\l\xi_1,\xi_2\r_{\Q} =\int_0^\infty\beta(s) \langle \xi_1(s),\xi_2(s)\rangle_ {H^2}ds.
\end{equation*}
Finally, we introduce the product space
\begin{equation*}
\H=H^{2}\times H^0\times H^0\times \M_\alpha\times \Q,
\end{equation*}
endowed with the norm
\begin{equation*}
\|(u,u_t,\vartheta,\eta,\xi)\|_{\H}^2=
\|u\|_{H^2}^2+\|u_t\|_{H^0}^2+\|\vartheta\|_{H^0}^2
+\|\eta\|_{\M_{\alpha}}^2+\|\xi\|_{\Q}^2.
\end{equation*}
In order to formulate the problem in the history space
setting we denote by $T$ and $T'$ the linear operators on $\M_\alpha$
and $\Q$ respectively, defined as
$$
T\eta=-\eta_s,\quad\eta\in{\D}(T), \qquad
T'\xi=-\xi_s,\quad\xi\in{\D}(T'),
$$
where $\D(T)=\{\eta\in\M_{\alpha}:\eta_s\in
\M_\alpha,\,\eta(0)=0\}$ and $\D(T')=\{\xi\in\Q:\xi_s\in
\Q,\,\xi(0)=0\}$, and $\eta_s$ (resp.\ $\xi_s$) stands for the
distributional derivative of $\eta$ (resp.\ $\xi$) with respect to
the internal variable $s$. On account of \eqref{K2A}, we immediately get
\begin{equation}
\label{etactrlA}
\l T\eta,\eta\r_{\M_\alpha}\leq 0, \qquad \l T'\xi,\xi\r_{\Q}\le 0,
\end{equation}
for $\eta\in\D(T)$ and $\xi\in\D(T')$. Let us introduce the formulation of
the problems. On account of the notation introduced above,
given $(u_0,u_1,\vartheta_0,\eta_0,\xi_0)$ in $\H$, find
$(u,u_t,\vartheta,\eta,\xi)\in C([0,\infty),\H)$ solution to
\begin{equation}
\tag{${\mathcal P}_{\alpha,\sigma}$}
\label{eq-mobileA}
\begin{cases}
\displaystyle
u_{tt}+\int_0^\infty \beta(s)A^2\xi(s)ds+A^2u-A^\sigma\vartheta=0, \\
\noalign{\vskip5pt} \displaystyle
\vartheta_t+\int_0^\infty \mu(s)A^\alpha\eta(s)ds+A^\sigma u_t=0, \\
\noalign{\vskip6pt}
\eta_t=T\eta+\vartheta, \\
\noalign{\vskip5pt} \xi_t=T'\xi+u_t,
\end{cases}
\end{equation}
for $t\in\R^+$, with initial conditions
$$
(u(0),u_t(0),\vartheta(0),\eta^0,\xi^0)=(u_0,u_1,\vartheta_0,\eta_0,\xi_0),
$$
and abstract boundary conditions
$$
u(t)\in\D(A^2),\qquad\vartheta(t)\in\D(A^\alpha),\qquad t\geq 0.
$$
System ${\mathcal P}_{\alpha,\sigma}$ allows us to provide a
description of the solutions in terms of a strongly continuous
semigroup of operators on $\H$. Indeed, setting
$\zeta(t)=(u(t),v(t),\vartheta(t),\eta^t,\xi^t)^{\top}$, the
problem rewrites as
\begin{equation*}
\frac{d}{dt}\zeta={\mathcal L}\zeta, \qquad
\zeta(0)=\zeta_0,
\end{equation*}
where ${\mathcal L}$ is the linear operator defined by
\begin{equation}
\label{cauchyA}
{\mathcal L}
\begin{pmatrix}
u \\
v \\
\vartheta \\
\eta \\
\xi
\end{pmatrix}
=
\begin{pmatrix}
v \\
-\int_0^\infty \beta(s)A^2\xi(s)ds-A^2u+A^\sigma\vartheta \\
\noalign{\vskip4pt}
-\int_0^\infty\mu(s)A^\alpha\eta(s)ds-A^\sigma v \\
\noalign{\vskip4pt}
\vartheta +T\eta \\
\noalign{\vskip4pt} v +T'\xi
\end{pmatrix}
\end{equation}
with domain
$$
{\mathcal D}({\mathcal L})= \left\{z\in\H \left|
\begin{aligned}
& Au\in H^2,\,\,
v\in H^{2\sigma},\,\,
\vartheta \in H^{2\sigma}  \\
& \textstyle\int_0^\infty\mu(s)A^\alpha \eta(s)ds\in H^0 \\
& \textstyle\int_0^\infty\beta(s)A^2(s)ds\in H^0 \\
& \eta\in{\D}(T),\,\,\, \xi\in{\D}(T')
\end{aligned}
\right\}. \right.
$$
\vskip4pt
\noindent
Since, by \eqref{etactrlA}, ${\mathcal L}$ is a
dissipative operator, arguing, e.g.,\ as in \cite{giopat}, and assuming
that \eqref{K0A}-\eqref{K2A} hold, we learn that ${\mathcal P}_{\alpha,\sigma}$
induces a $C_0$-semigroup $S_{\alpha,\sigma}(t)$ of contractions on $\H$.
\vskip4pt
\noindent
Assuming that \eqref{K0A}-\eqref{add-int} hold, we have the following

\begin{theoremA}[{\bf pointwise decay}]
\label{dipp}
For every $\alpha,\sigma\geq 0$,
$$
\lim_{t\to\infty} \|S_{\alpha,\sigma}(t)z_0\|_{\H}=0,
\qquad\forall\, z_0\in\H.
$$
\end{theoremA}

\noindent
The main result of the Appendix are the following

\begin{theoremA}[{\bf non-exponential decay I}]
\label{ThmMAIN0}
Assume that $0\leq\alpha<2$, $\sigma\geq 0$ and
\begin{equation*}
\mu(s)=\kappa_1s^{-\omega_1}e^{-\delta_1s},\quad
0\leq \omega_1<\frac{2-\alpha}{2},\,\,\,\kappa_1,\delta_1>0,
\qquad\,\,\,\, \beta(s)=0.
\end{equation*}
Then $S_{\alpha,\sigma}(t)$ is not exponentially stable on $\H$.
\end{theoremA}

\begin{theoremA}[{\bf non-exponential decay II}]
\label{ThmMAIN}
Assume that for some $\kappa_1,\kappa_2,\delta_1,\delta_2>0$,
\begin{equation*}
\mu(s)=\kappa_1 s^{-\omega_1}e^{-\delta_1s}, \qquad
\beta(s)=\kappa_2 s^{-\omega_2}e^{-\delta_2 s}.
\end{equation*}
Furthermore, suppose that
$$
0\leq\alpha<2,\qquad 0\leq\sigma<1, \qquad \alpha\leq 2\sigma,
$$
and that
$$
\frac{2\sigma-\alpha}{2}\leq \omega_1<\frac{2-\alpha}{2},
\qquad
0\leq \omega_2\leq\omega_1-\frac{2\sigma-\alpha}{2}.
$$
\vskip1pt
\noindent
Then $S_{\alpha,\sigma}(t)$ is not exponentially stable on $\H$.
\end{theoremA}

\begin{remarkA}
The condition $\alpha<2\sigma$ on the $A$ powers is
not new in thermoelasticity.
It appears for instance (among other restrictions) in the study of
smoothing/non-smoothing properties for a class of abstract (memory free)
thermoelastic systems (cf.\ \cite{riv-rak2}).
\end{remarkA}

\begin{remarkA}
The memory kernels of Theorems \ref{ThmMAIN0}-\ref{ThmMAIN}
also satisfy the extra summability condition \eqref{add-int}.
Hence, any trajectory of the system goes to zero, but with
an arbitrarily slow decay rate, according to the chosen initial data.
\end{remarkA}

The rest of the Appendix is devoted to the proof of the above results.

\subsection{Proof of Theorem \ref{dipp}}
In order to prove the result, we shall exploit the following sufficient
condition for the (pointwise) decay to zero of any trajectory of a
linear gradient system (cf., e.g., \cite[Thm. A.2 and Cor. A.3]{GP-Riv}).

\begin{lemmaA}
\label{absdeca}
Let $S(t)$ be a linear gradient system on a Banach
space $\H$, let $z_0\in\H$ and assume that
$$
\text{$\bigcup_{t\geq 0}S(t)z_0$ is relatively compact in $\H$}.
$$
Then
$$
\lim_{t\to\infty} S(t)z_0=0.
$$
The same holds if the hypotheses are satisfied for all $z_0\in\X$,
with $\X$ dense subset of $\H$.
\end{lemmaA}
\vskip2pt
By exploiting \eqref{K0A}-\eqref{K2A} and \eqref{etactrlA} and observing that, by
\eqref{add-int}, $\mu$ and $\beta$ cannot be identically equal to
zero, it is easily seen that $S_{\alpha,\sigma}(t)$ is a gradient
system on $\H$ (argue, e.g.,\ as in~\cite[Prop. 3.2]{GP-Riv}).
We shall also set
$$
\M_\alpha^1=L^2_{\mu}(\R^+,H^{\alpha+1}), \quad
\Q^1=L^2_{\beta}(\R^+,H^3), \quad \H^1=H^{3}\times H^{1}\times
H^{1}\times \M_\alpha^1\times \Q^1.
$$
If ${\mathcal L}$ denotes the linear operator defined in
\eqref{cauchyA}, since the space ${\mathcal D}({\mathcal
L})\cap\H^1$ is dense in $\H$, according to Lemma~\ref{absdeca} it
is sufficient to check the assumptions for a fixed $z_0\in
{\mathcal D}({\mathcal L})\cap\H^1$. Let $C=C(z_0)$ denote a
generic positive constant. It is readily seen that
$\|S_{\alpha,\sigma}(t)z_0\|_{\H^1}\leq C$ for all $t\geq 0$.
Indeed, by \eqref{etactrlA} it suffices to
multiply the equations of ${\mathcal P}_{\alpha,\sigma}$ by $u_t$
in $H^{1}$, by $\vartheta$ in $H^{1}$, by $\eta$ in
$\M^1_\alpha$ and by $\xi$ in $\Q^1$ respectively and add the
resulting equations. Let us consider the sets
$$
\CC_1=\overline{\bigcup_{t\geq 1}\eta^t}^{\M_\alpha} \qquad
\text{and} \qquad \CC_2=\overline{\bigcup_{t\geq 1}\xi^t}^{\Q}.
$$
We claim that $\CC_1\times\CC_2\subset\M_\alpha^1\times\Q^1$ is
compactly embedded into $\M_\alpha\times\Q$. To this aim, we
recall the following compactness result (see, e.g., \cite[Lemma
2.1]{GP-Riv}) for the spaces $\M^1_\alpha\times\Q^1$. Assume that
$\CC_1\subset\M^1_\alpha$ and $\CC_2\subset\Q^1$ satisfy:
\vskip5pt
\begin{itemize}
\item[{\rm (i)}] $\sup\limits_{\eta\in\CC_1} \|\eta\|_{\M_\alpha^1}<\infty$
\qquad \text{and} \qquad
$\sup\limits_{\eta\in\CC_1}\|\eta_s\|_{\M_\alpha}<\infty$,
\vskip6pt
\item[{\rm (ii)}] $\sup\limits_{\xi\in\CC_2} \|\xi\|_{\Q^1}<\infty$
\qquad \text{and} \qquad
$\sup\limits_{\xi\in\CC_2}\|\xi_s\|_{\Q}<\infty$, \vskip6pt
\item[{\rm (iii)}] $\displaystyle\lim_{x\to\infty}\big[
\sup\limits_{\eta\in\CC_1}\T_\eta(x)\big]=0$ \qquad \text{and}
\qquad $\displaystyle\lim_{x\to\infty}\big[
\sup_{\xi\in\CC_2}\T_\xi(x)\big]=0$,
\end{itemize}
where the tails functions $\T_\eta$ and $\T_\xi$ are defined by
\begin{align*}
\T_\eta(x)
&=\int_{(0,1/x)\cup(x,\infty)}\mu(s)\|A^{\alpha/2}\eta(s)\|^2 ds,
\qquad x\geq 1, \\
\noalign{\vskip3pt}
\T_\xi(x)
&=\int_{(0,1/x)\cup(x,\infty)}\beta(s)\|A\xi(s)\|^2 ds, \qquad x\geq 1.
\end{align*}
Then $\CC_1\times\CC_2$ is relatively compact in
$\M_\alpha\times\Q$. Indeed, by simply mimicking the proofs of
\cite[Lemma 4.3 and Lemma 4.4]{GP-Riv}, exploiting the
representation formulas for $\eta^t$ and $\xi^t$
\begin{align*}
\eta^t(s)&=
\left\{\begin{array}{ll}
\int_0^s \th(t-y)dy, & \quad \mbox{$0<s\le t$},\\
\noalign{\vskip6pt}
\eta_0(s-t)+\int_0^t \th(t-y)dy, & \quad \mbox{$s>t$,}
\end{array}\right. \\
\noalign{\vskip5pt}
\xi^t(s)&=\left\{\begin{array}{ll}
u(t)-u(t-s), & \quad \mbox{$0<s\le t$},\\
\noalign{\vskip6pt}
\xi_0(s-t)+u(t)-u(0), & \quad \mbox{$s>t$,}
\end{array}\right.
\end{align*}
it is readily seen that (i)-(iii) are fulfilled (we point out
that the addition summability assumption \eqref{add-int} on $\mu$
pops up in the proof of (iii) for $\eta^t$). Now, consider the set
$$
{\mathcal K}=B_{H^3\times H^{1}\times
H^{1}}(C)\times\CC_1\times\CC_2.
$$
Then, ${\mathcal K}$ is compact in $\H$ being
$B_{H^3\times H^{1}\times H^{1}}(C)$ compact in
$H^2\times H^0\times H^0$ and $\CC_1\times\CC_2$ compact in
$\M_\alpha\times\Q$. Moreover, by construction, there holds
$S_{\alpha,\sigma}(t)z_0\in{\mathcal K}$ for every
$t\geq 0$. Therefore, by Lemma~\ref{absdeca}, we have
$S_{\alpha,\sigma}(t)z_0\to 0$ in $\H$ as $t\to\infty$. \qed

\subsection{Proof of Theorems \ref{ThmMAIN0} and \ref{ThmMAIN}}

To prove the results, we shall exploit the following
classical result due to Pr\"uss \cite{PRU}.

\begin{lemmaA}
\label{lemmapr}
Let $S(t)=e^{t{\mathcal L}}$ be a $C_0$-semigroup of contractions
on a  Hilbert space $\H$. Then $S(t)$ is exponentially stable if and only if
$i\R$ belongs to the resolvent set of ${\mathcal L}$, and
there exists $\eps>0$ such that
\begin{equation*}
\inf_{\lambda\in\R}\|(i\lambda\I-{\mathcal L})z\|_{\H} \geq
\eps\|z\|_{\H},\qquad\forall z\in\D({\mathcal L}).
\end{equation*}
\end{lemmaA}
\vskip2pt
We start with the proof of Theorem \ref{ThmMAIN}, for the proof of Theorem \ref{ThmMAIN0}
is just a simple by-product. Let ${\mathcal L}$ be the linear operator defined in
\eqref{cauchyA}. For $\lambda\in\R$ and for $\tilde z=(0,0,0,\teta,\tilde\xi)^\top\in\H$,
we consider the complex equation $(i\lambda\I-{\mathcal L})z=\tilde z$, which explicitly
writes as
$$
\begin{cases}
 i\lambda u - v=0,  &\\
\noalign{\vskip4pt}
 i\lambda v+\displaystyle\int_0^\infty \beta(s) A^2\xi(s) ds+A^2u-A^\sigma\th=0,  &\\
\noalign{\vskip4pt}
 i\lambda\th+\displaystyle\int_0^\infty \mu(s) A^\alpha\eta(s)ds+A^\sigma v=0,   &\\
\noalign{\vskip4pt}
 i\lambda\xi -v +\xi_s =\tilde\xi, \\
\noalign{\vskip6pt}
 i\lambda\eta -\th +\eta_s =\teta.
\end{cases}
$$
We shall denote by $\{\gamma_n\}$ the sequence of (positive) eigenvalues
of $A$ and by $\{w_n\}$ the corresponding sequence of normalized
eigenvectors. We choose
$$
\teta(s)=\teta_n(s)= \gamma_n^{-\alpha/2}w_n, \qquad
\tilde\xi(s)=\tilde\xi_n(s)= \Lambda_n w_n.
$$
where $\Lambda_n=\Lambda_n(\gamma_n)$ will be suitably chosen later on. If
$\tilde z_n=(0,0,0,\teta_n,\tilde\xi_n)^\top$, then it holds
\begin{equation}
\label{initcd}
\|\tilde z_n\|_{\H}^2=k_0+h_0(\gamma_n\Lambda_n)^2,
\qquad\text{for all $n\in\N$},
\end{equation}
where we have set
\begin{equation*}
k_0 =\int_0^\infty\mu(s)ds, \qquad h_0 =\int_0^\infty\beta(s)ds.
\end{equation*}
We shall prove the assertion by applying Lemma \ref{lemmapr},
arguing by contradiction. To this aim, we find a
sequence $\{\lambda_n\}$ in $\R$ and a corresponding solution
$z_n$ such that $\|z_n\|_{\H}\to\infty$, as $n\to\infty$. We
search for a solution $z=(u,v,\th,\eta,\xi)^\top$ of the form
$$
u=u_n=pw_n,\,\,\,
v=v_n=qw_n,\,\,\,
\th=\th_n=rw_n,\,\,\,
\eta=\eta_n=\varphi w_n,\,\,\,
\xi=\xi_n=\psi w_n,
$$
where $p,q,r\in\C$, $\varphi\in H^1_\mu(\R^+)$ and $\psi \in H^1_\beta(\R^+)$, with
$\varphi(0)=\psi(0)=0$. Whence, the above system leads to the
following equations
$$
\begin{cases}
 i\lambda p-q=0,  &\\
 \gamma_n^2 p-\lambda^2 p-\gamma_n^\sigma r+\gamma_n^2
\displaystyle\int_0^\infty\beta(s) \psi(s) ds=0, &\\
\noalign{\vskip5pt}
 i\lambda r+i\lambda\gamma_n^\sigma p
+\gamma_n^\alpha\displaystyle\int_0^\infty\mu(s) \varphi(s) ds=0, &\\
\noalign{\vskip5pt}
 i \lambda \varphi(s)-r+\varphi_s(s)
=\frac{1}{\gamma_n^{\alpha/2}}, &\\
\noalign{\vskip5pt}
 i \lambda \psi(s)-q+\psi_s(s)=\Lambda_n.
\end{cases}
$$
Imposing $\varphi(0)=\psi(0)=0$, we can integrate the last two
equations, getting
\begin{align*}
\varphi(s)&=\frac{1}{i\lambda}
\big(r+\gamma_n^{-\alpha/2}\big)\big(1-e^{-i\lambda s}\big), \\
\noalign{\vskip3pt}
\psi(s)&=\frac{1}{i\lambda}(q+\Lambda_n)\big(1-e^{-i\lambda
s}\big).
\end{align*}
Then, we are led to the following system
\begin{equation}
\label{sist-compl}
\begin{cases}
i\lambda r+i\lambda\gamma_n^\sigma p+
\frac{\gamma_n^\alpha}{i\lambda}\big(r+\gamma_n^{-\alpha/2}\big)(k_0-c(\lambda))=0, &\\
\noalign{\vskip6pt}
\gamma_n^2p-\lambda^2 p-\gamma_n^\sigma r+
\frac{\gamma_n^2}{i\lambda}(i\lambda p+\Lambda_n)(h_0-b(\lambda))=0,
\end{cases}
\end{equation}
being $c(\lambda)$ and $b(\lambda)$ the Laplace transform
of the kernels $\mu$ and $\beta$ respectively,
\begin{align*}
c(\lambda) &=\int_0^\infty\mu(s) e^{-i\lambda s}ds,\qquad\lambda\in\R^+ \\
\noalign{\vskip3pt}
b(\lambda) &=\int_0^\infty\beta(s) e^{-i\lambda s}ds,\qquad\lambda\in\R^+.
\end{align*}
We now impose the conditions
\begin{equation}
\label{impose}
i\lambda r+i\lambda\gamma_n^\sigma p+
\frac{rk_0\gamma_n^\alpha}{i\lambda}=0, \qquad\,
p=\frac{\gamma_n^\sigma r}{(1+h_0)\gamma_n^2-\lambda^2}.
\end{equation}
These yield the fourth order algebraic equation
\begin{equation}
\label{fourto}
\lambda^4-\big[(1+h_0)\gamma_n^2+\gamma_n^{2\sigma}+k_0\gamma_n^\alpha\big]\lambda^2
+k_0(1+h_0)\gamma^{\alpha+2}_n=0.
\end{equation}
Taking into account that, since $\alpha<2$ and $\sigma<1$, we have
$$
(1+h_0)\gamma_n^2+\gamma_n^{2\sigma}+k_0\gamma_n^\alpha={\mathcal O}(\gamma_n^2),\qquad\text{as
$n\to\infty$},
$$
it is easy to realize that \eqref{fourto} admits a real positive solution
$\lambda=\lambda_n=\lambda_n(\gamma_n)={\mathcal O}(\gamma_n)$, as $n\to\infty$.
Consequently, setting $c_n=c(\lambda_n)$ and $b_n=b(\lambda_n)$, from
\eqref{sist-compl}-\eqref{impose}, we get
\begin{align*}
r&=\frac{k_0-c_n}{\gamma_n^{\alpha/2}c_n}=r_n(\gamma_n),  \\
\noalign{\vskip2pt}
p&=\frac{\gamma_n^\sigma r_n}{(1+h_0)\gamma_n^2-\lambda_n^2}=p_n(\gamma_n), \\
\noalign{\vskip6pt}
\Lambda_n&=\frac{i\lambda_np_nb_n}{h_0-b_n}=\Lambda_n(\gamma_n).
\end{align*}
Notice that the above quantities depend solely
on the eigenvalues $\gamma_n$ of $A$. Moreover,
\begin{align*}
b_n=\int_0^\infty
\kappa_2s^{-\omega_2}e^{-(i\lambda_n+\delta_2)s}ds
&=\kappa_2\lambda_n^{\omega_2-1}\Big(i+\frac{\delta_2}{\lambda_n}\Big)^{\omega_2-1}
\Gamma(1-\omega_2) \\
\noalign{\vskip2pt}
& ={\mathcal O}\left(\lambda_n^{\omega_2-1}\right)
={\mathcal O}\left(\gamma_n^{\omega_2-1}\right), \\
\noalign{\vskip5pt}
c_n=\int_0^\infty \kappa_1s^{-\omega_1}e^{-(i\lambda_n+\delta_1)s}ds&=
\kappa_1\lambda_n^{\omega_1-1}\Big(i+\frac{\delta_1}{\lambda_n}
\Big)^{\omega_1-1}\Gamma(1-\omega_1) \\
\noalign{\vskip2pt}
& ={\mathcal O}\left(\lambda_n^{\omega_1-1}\right)={\mathcal O}
\left(\gamma_n^{\omega_1-1}\right),
\end{align*}
as $n\to\infty$, where $\Gamma$ is the Gamma Function, so that
\begin{equation*}
\frac{b_n}{c_n}
={\mathcal O}\left(\gamma_n^{\omega_2-\omega_1}\right)\qquad \text{as
$n\to\infty$}.
\end{equation*}
As a consequence, we obtain
\begin{align}
\label{Ocap}
\gamma_n\Lambda_n(\gamma_n)&=\frac{i\lambda_nb_n}{h_0-b_n}\frac{\gamma_n^{\sigma+1}}{
(1+h_0)\gamma_n^2-\lambda_n^2}\frac{k_0-c_n}{\gamma_n^{\alpha/2}c_n}  \\
\noalign{\vskip4pt}
&={\mathcal O}\bigg(\frac{\gamma_n^{\sigma+2-\alpha/2}}{
(1+h_0)\gamma_n^2-\lambda_n^2}\frac{b_n}{c_n}\bigg)
={\mathcal O}\big(\gamma_n^{\omega_2-\omega_1+\sigma-\alpha/2}\big), \notag
\end{align}
as $n\to\infty$. By \eqref{initcd}, \eqref{Ocap} and the assumptions
on $\sig,\alpha,\omega_1,\omega_2$, we learn that
$$
\sup_{n\geq 1}\|\tilde z_n\|_{\H}<\infty.
$$
On the other hand, by the assumptions on $\sig,\alpha,\omega_1$,
we have $|r_n|\to \infty$ as $n\to\infty$, yielding
$$
\|z_n\|_{\H}\geq \|\th_n\|=|r_n|\to\infty,\qquad\text{as $n\to\infty$},
$$
which readily yields a contradiction
and concludes the proof of Theorem \ref{ThmMAIN}.
\vskip2pt
\noindent
The proof of Theorem \ref{ThmMAIN0} simply follows by mimicking
the above steps, observing that by assumption we have $h_0=0$.
In particular $\|\tilde z_n\|_{\H}=\sqrt{k_0}$ by \eqref{initcd}, whereas
$\omega_1<\frac{2-\alpha}{2}$ implies that $\|z_n\|_{\H}\to\infty$
as $n\to\infty$, yielding again the assertion. \qed

\vskip30pt

\noindent
{\bf Acknowledgment.} The authors are grateful to the referees for their appropriate and
useful remarks.

\vskip40pt

\vskip50pt

\begin{thebibliography}{99}

\bibitem{cps1}
{\au M.~Conti, V.~Pata, M.~Squassina},
{\ti Singular limit of differential systems with memory},
{\jou Indiana Univ.\ Math.\ J.}
\no{55}{169--216}{2006}

\bibitem{cps2}
{\au M.~Conti, V.~Pata, M.~Squassina},
{\ti Singular limit of dissipative hyperbolic equations with memory},
{\jou Discrete Contin.\ Dyn.\ Syst.}
\no{{\rm Special Issue}}{200--208}{2005}



\bibitem{gnp}
{\au C.~Giorgi, M.G.~Naso, V.~Pata},
{\ti Energy decay of electromagnetic systems with memory},
{\jou Math.\ Models Methods Appl.\ Sci.}
\no{15}{1489--1502}{2005}

\bibitem{giopat}
{\au C.~Giorgi, V.~Pata},
{\ti Stability of linear thermoelastic systems with memory},
{\jou Math.\ Models Methods Appl.\ Sci.}
\no{11}{627--644}{2001}

\bibitem{GP-Riv}
{\au M.~Grasselli, J.E.~Mu\~noz Rivera, V.~Pata},
{\ti On the energy decay of the linear thermoelastic plate with memory},
{\jou J.\ Math.\ Anal.\ Appl.}
\no{309}{1--14}{2005}

\bibitem{Terreni}
{\au M.~Grasselli, V.~Pata},
{\ti Uniform attractors of nonautonomous systems with memory},
in ``Evolution Equations, Semigroups and Functional Analysis''
(A.~Lorenzi and B.~Ruf, Eds.),
\eds{pp.155--178, Progr.\ Nonlinear Differential Equations
Appl.\ no.50, Birkh\"{a}user}{Boston}{2002}

\bibitem{GS1}
{\au M.~Grasselli, M.~Squassina},
{\ti Exponential stability and singular limit for a linear
thermoelastic plate with memory effects},
{\jou Adv.\ Math.\ Sci.\ Appl.}
\no{16}{15--31}{2006}

\bibitem{GuPi}
{\au M.E.~Gurtin, A.C.~Pipkin},
{\ti A general theory of heat conduction with finite wave speeds},
{\jou Arch.\ Rational Mech.\ Anal.}
\no{31}{113--126}{1968}

\bibitem{hazu}
{\au A.~Haraux, E.~Zuazua},
{\ti Decay estimate for some damped hyperbolic equations},
{\jou Arch.\ Rational Mech.\ Anal.}
\no{100}{191--208}{1988}

\bibitem{LaTr}
{\au I.~Lasiecka and R.~Triggiani},
{\ti Control theory for partial differential equations: continuous and approximation theories. I. 
Abstract parabolic systems}, 
\eds{Cambridge University Press}{Cambridge}{2000}








\bibitem{riv-rak}
{\au J.E.~Mu\~noz Rivera, R.~Racke},
{\ti Smoothing properties, decay and global existence of solutions
to nonlinear coupled systems of thermoelastic type},
{\jou SIAM J.\ Math.\ Anal.}
\no{26}{1547--1563}{1995}

\bibitem{riv-rak2}
{\au J.E.~Mu\~noz Rivera, R.~Racke},
{\ti Large solution and smoothing properties for nonlinear thermoelastic systems},
{\jou J.\ Differential Equations}
\no{127}{454--483}{1996}


\bibitem{naso}
{\au M.G.~Naso},
{\ti Exponential stability of a viscoelastic plate with thermal memory},
{\jou Riv.\ Mat.\ Univ.\ Parma}
\no{3}{37--56}{2000}

\bibitem{PRU}
{\au J.~Pr\"uss},
{\ti On the spectrum of $C_0$-semigroups},
{\jou Trans.\ Amer.\ Math.\ Soc.}
\no{284}{847--857}{1984}


\bibitem{zua}
{\au E.~Zuazua},
{\ti Stability and decay for a class of nonlinear hyperbolic problems},
{\jou Asymptotic Anal.}
\no{1}{1--28}{1998}

\end{thebibliography}
\end{document}